\documentclass[review]{elsarticle}
\usepackage{lineno,hyperref}
\modulolinenumbers[5]
\usepackage{flushend}
\usepackage[noend]{algorithmic} 
\usepackage{algorithm,caption}
\algsetup{indent=2em} 
 
\usepackage{subfig}
\restylefloat{algorithm}
\usepackage{longtable}
\usepackage{multirow}
\usepackage{lscape}
\usepackage{amsmath}
\usepackage{amssymb}
\usepackage{graphicx}
\usepackage{xcolor}
\usepackage{amsfonts}
\usepackage[utf8]{inputenc}
\newtheorem{Proof}{Proof}[section]

\newtheorem{Remark}{Remark}[section]

\newtheorem{definition}{Definition}[section]

\newtheorem{Theorem}{Theorem}[section]
\journal{Journal of \LaTeX\ Templates}

\usepackage{setspace,lipsum}

\begin{document}

\begin{frontmatter}

\title{ Spatial Metric Space for Pattern Recognition Problems}
\tnotetext[mytitlenote]
{Corresponding Author: Q.M. Danish Lohani\\
 ~~~~~Email: danishlohani@cs.sau.ac.in  
	}

\author{$^{1}$Q.M. Danish Lohani, $^{2}$Ashutosh Tiwari, $^{3}$Mohd Shoaib Khan}
\address{$^{1,2}$Department of Mathematics, South Asian University, New Delhi 110068, India\\
$^{3}$ Department of Data Sciences and Analytics, Ramaiah
University of Applied Sciences, Bangalore, India\\
Email: $^{1}$danishlohani@cs.sau.ac.in,$^{2}$ashutoshtiwari591@gmail.com,$^{3}$shoaibkhanamu@gmail.com}

\begin{abstract}
 The definition of weighted distance measure involves weights. The paper proposes a weighted distance measure without the help of weights. Here, weights are intrinsically added to the measure, and for this, the concept of metric space is generalized based on a novel divided difference operator. The proposed operator is used over a two-dimensional sequence of bounded variation, and it generalizes metric space with the introduction of a multivalued metric space called spatial metric space. The environment considered for the study is a two-dimensional Atanassov intuitionistic fuzzy set (AIFS) under the assumption that membership and non-membership components are its independent variables. The weighted distance measure is proposed as a spatial distance measure in the spatial metric space. The spatial distance measure consists of three branches. In the first branch, there is a domination of membership values, non-membership values dominate the second branch, and the third branch is equidominant. The domination of membership and non-membership values are not in the form of weights in the proposed spatial distance measure, and hence it is a measure independent of weights. The proposed spatial metric space is mathematically studied, and as an implication, the spatial similarity measure is multivalued in nature. The spatial similarity measure can recognize a maximum of three patterns simultaneously. The spatial similarity measure is tested for the pattern recognition problems and the obtained classification results are compared with some other existing similarity measures to show its potency. This study connects the double sequence to the application domain via a divided difference operator for the first time while proposing a novel divided difference operator-based spatial metric space.
\end{abstract}

\begin{keyword}
\texttt{}Atanassov intuitionistic fuzzy sets, Bounded variation, Similarity measure, Pattern recognition.
\end{keyword}

\end{frontmatter}

\section{Introduction} \label{sec1}
The difference operator helps in defining a well-known sequence of bounded variation. The divided difference operator is popularly used over a double sequence. 
The double sequence is an extension of a sequence containing two indices and reduces to a sequence if its one index is kept constant. Sometimes double sequence is not partitioned, which yields an infinite ordered matrix. In order to structure the double sequence, a functional relationship of two indices is used, such that the first index performs the numbering of partitions, and the second index informs about the total elements present in each partition. The convergence of a double sequence is not uniquely defined. Pringsheim \cite{40} proposed the convergence of double sequence for the first time. In sequence, convergence implies boundedness, which does not hold for the double sequence. The double sequence fails to match a single sequence in the application domain because double sequence-based divided difference operators are not utilized enough in real-world applications. So, the concept of double sequence could not attract the researchers of engineering and applied sciences, although it differs theoretically from the sequence. The double sequences are also studied in the fuzzy environment (see: \cite{17}, \cite{41}, \cite{49}, \cite{50}, \cite{37}, \cite{31}, \cite{32}, \cite{33}).

The difference operator is used as a tool to estimate the patterns of arc length, and it involves the concept of bounded variation. There are several theoretical studies about bounded variation in sequence spaces, but it is rarely used for pattern recognition \cite{2},\cite{26}. The fuzzy bounded variation motivates for the pattern analysis in the two-dimensional plane. Atanassov intuitionistic fuzzy set (AIFS) \cite{5,4} resolves the uncertainty arising due to impreciseness and incompleteness. AIFS is considered to be an element of the two-dimensional plane in some studies of decision-making, image processing, and pattern recognition (see: \cite{3},\cite{47},\cite{54},\cite{2},\cite{9},\cite{56}). 
The AIFS based research work of image processing, decision-making, and pattern recognition involves a similarity measure. Similarity measure, basically is a tool to compute the similarity between the elements. In the AIFS domain, similarity and distance measures mostly complement each other \cite{19}, and this result motivates us to pursue an AIFS based study. Here, a distance based similarity measure is introduced for two-dimensional AIFS. The $p$-summable or $p$-summable bounded variation based AIFS distance measures are highly used as they have simple mathematical formulations. Usually, Hamming and Euclidean distances \cite{5,8} are exploited to address the problems of pattern recognition as the geometrical structure significantly differs from each other. The AIFS based $p$-summable distance measure could not be differentiated from its double sequence based counterpart due to the equal computation of the distances between any two elements \cite{2} and hence the double sequence based counterpart does not help in the pattern recognition problem. The distance measures defined over bounded variation sequence involve a difference operator for the computation of the component-wise absolute difference, and its double sequence counterpart uses a divided difference operator to compute the same differences. We have utilized a divided difference operator in a two-dimensional AIFS plane for pattern recognition while proposing an AIFS similarity measure. Some of the existing well-known AIFS similarity measures are given in Table-\ref{similarity}. There are some problems in which the three-dimensional aspect of AIFS is explored (see: \cite{46},\cite{42},\cite{48},\cite{45}).  
 
 The aim of the paper is to substantiate multivalued analysis under a metrizable space. The metrizibility in the absolutely summable bounded variation involves a difference operator. The difference operator operates over the scaler/real valued sequences. In real valued double sequences, the divided difference operator deals with the two independent real variables. The divided difference operator is modified in contemplation to two-dimensional real valued sequences. The concept of bounded variation is exploited to propose multivalued absolutely summable bounded variation. The proposed multivalued bounded variation gives rise to a concept of multivalued metric space. The proposed multivalued metric space is three branched and hence it is named as spatial metric space. None of the branches of the spatial metric space involves weights and hence spatial distance measure is independent of weights. We summarize this discussion while saying that a spatial metric is a tool of multivalued analysis. The convergence is proved to verify the multivalued analysis in spatial metric space. Here, our major concern is to establish the multivalued analysis mathematically, and beside this some pattern recognition applications are also provided.

 Now, we summarize the major offerings of the paper as:
 \begin{itemize}
	\item[(1)] A geometrical interpretation of the difference operator is given in terms of one point-based neighbourhood, and involvement of three point-based neighbourhood is observed in the divided difference operator over a real-valued set.
	\item[(2)] The selection of a three-point based neighbourhood is easy in the two-dimensional real-valued space, so we have considered a two-dimensional space called the AIFS plane. Membership and non-membership components are taken as the two dimensions of the AIFS.
	\item[(3)] A corroboration of the divided difference operator and bounded variation in the AIFS plane resulted in the definition of multivalued absolutely summable bounded variation.
	\item[(4)] In multivalued absolutely summable bounded variation, a spatial metric space of bounded variation has been proposed along with its convergence. 
	\item[(5)] The strength of the novel divided difference based spatial metric space lies in the intrinsic assignment of weights, and it has been highlighted over pattern recognition problems due to a proposed spatial similarity measure (SSM). 
	\item[(6)] The distance measure of the spatial metric space is used to propose a three branched SSM. The three branches of SSM are termed as membership dominant, non-membership dominant, and equidominant. The SSM is independent of weights. Moreover, the proposed SSM has successfully resolved the cancer diagnosis and medical diagnosis problems.
	\item[(7)] This study is essential as it negates the application lagging of the double sequence. Here, we have used a novel divided difference operator of the double sequence to establish the connection between the fuzzy and real domains.
\end{itemize}

The paper is divided into six sections. Section \ref{sec1} contains the introduction of the research work. In Section \ref{sec2}, definitions along with the notations used throughout the paper are mentioned. Section \ref{sec3} presents the proposal of a spatial metric space, a novel spatial similarity measure, and some mathematical properties related to them. In Section \ref{sec4}, the proposed spatial similarity measure is applied to address well-known pattern recognition problems. The conclusion and future study are in are in Section \ref{sec5} and Section \ref{sec6} respectively.

\begin{table}[htp]
	\renewcommand{\arraystretch}{0.5}
	\scriptsize
	\centering
	\caption{Some well-known similarity measures for AIFSs with notation $(\mu_{B}(y_j),\nu_{B}(y_j))=(a_j,b_j)$ and $(\mu_{A}(y_j),\nu_{A}(y_j))=(c_j,d_j)$ }
	\label{similarity}
	\begin{tabular}{l|l}
		
		\hline
		
		~\\
		
		~\\
		
		\hspace{0.5cm}  Source  & \hspace{5cm} Similarity measures\\
		
		~\\
		
		\hline
		
		Atanassov~\cite{5} & $S_{ \rm L_P}=1-\sqrt[p]{\frac{1}{2k}\sum\limits_{j=1}^{k}(a_j-c_j)^p+(b_j-d_j)^p}$ where, $1\leq p<\infty$\\
		
		\hline
		
		Boran et al. \cite{8} &$S_{\rm Bd}=1-\sqrt[p]{\frac{1}{2k(z+1)^p}\sum\limits_{j=1}^{k}(| z(a_j-c_j)-(b_j-d_j)|^p+ | z(b_j-d_j)-(a_j-c_j)|^p}~$\\& where, $1\leq p<\infty$ and $2\leq z<\infty$ \\
		
		\hline
		
		Chen et al. \cite{10} & $S_{\rm C}=1-\frac{1}{2k}\sum\limits_{j=1}^{k}|(a_j-b_j)-(c_j-d_j)|$\\
		
		\\
		
		\hline
		

		
		
		
		
		Dengfeng et al. \cite{14}  & $S_{\rm Dc}=1-\sqrt[p]{\frac{1}{k}\sum\limits_{j=1}^{k}\left |  m_{B}(y_j)-m_{A}(y_j) \right |^p}, ~where,~ m_B(y_j)=\frac{a_j+1-b_j}{2},~ m_A(y_j)=\frac{c_j+1-d_j}{2}$\\
		
		\hline
		
		Fan et al. \cite{15} &  $S_{\rm Fz}=1-\frac{1}{4k}(\sum\limits_{j=1}^{k}|(a_j-b_j)-(c_j-d_j)|+|(a_j-c_j)-(b_j-d_j)|)$\\
		
		\hline
		
		
		
		Hong et al. \cite{18} & $S_{\rm Hk}=1-\frac{1}{2k}(\sum\limits_{j=1}^{k}|a_j-c_j|+|b_j-d_j|)$\\
		
		\hline
		
		Hung et al. \cite{19} &$S_{\rm Hy_1}^1=1-\frac{1}{k}\sum\limits_{j=1}^{k}\max(|a_j-c_j|,|b_j-d_j|)$, $S_{\rm Hy_2}^1=\frac{e^{S_{Hy_1}-1}-e^{-1}}{1-e^{-1}}$, $S_{Hy_3}^1=\frac{S_{\rm Hy_1}}{2-S_{\rm Hy_1}}$\\
		
		\hline
		
		Hung et al. \cite{20}&$S_{Hy_1}^2=\frac{2^{1/p}-D_p}{2^{1/p}}$, $S_{\rm Hy_2}^2=\frac{e^{-D_p}-e^{-2^{1/p}}}{1-e^{-2^{1/p}}}$, and $S_{\rm Hy_3}^2=\frac{2^{1/p}-D_p}{2^{1/p}(1+D_p)}$; $D_p=\frac{1}{k}\sqrt[p]{\sum\limits_{j=1}^{k}|a_j-c_j|^p+|b_j-d_j|^p}$\\
		
		\hline
		
		Hung et al. \cite{21} &
		$S_{Hy}^3= \frac{1}{k}\sum\limits_{j=1}^{k}1-(\frac{1}{2}|a_j-c_j|+|b_j-d_j|)$\\ 

		\\
		
		\hline
		
		
		
		
		
		
		Li et al. \cite{27}    &$S_{\rm Lzd}=1-\sqrt{\frac{1}{2k}\sum\limits_{j=1}^{k}|a_j-c_j|^2+|b_j-d_j|^2}$\\
		
		\hline
		
		Liang et al. \cite{28}    &$S_{\rm Ls}=1-\sqrt[p]{\frac{1}{k}(\sum\limits_{j=1}^{k}|s_j-t_j|^p)}$, where $s_j=\frac{|a_j-c_j|}{2} $,~ $t_j=\frac{|(1-b_j)-(1-d_j)|}{2} $\\
		
		\hline
		
		

		Mitchell \cite{30} & $S_{\rm M}=\frac{1}{2}(\rho_{1}+\rho_{2})$;  $\rho_{1}=1-\sqrt[p]{\frac{1}{k}\sum\limits_{j=1}^{k}\left |a_j-c_j \right |^p}$, $\rho_{2}=1-\sqrt[p]{\frac{1}{k}\sum\limits_{j=1}^{k}\left |b_j-d_j \right |^p}$ \\
		
		\hline
		
		Nagan et al. \cite{34} & $S_{\rm Hm}=1-\frac{1}{3k}\sum\limits_{j=1}^{k}|a_j-c_j|+|b_j-d_j|+|\max(a_j,d_j)-\max(c_j,b_j)|$ \\
		
		\hline
		
		Ashraf et al. \cite{2} & $S_{\rm Az_{p}}= 1-\frac{1}{\sqrt[p]{2k}}(|a_1-c_1|^p+|b_1-d_1|^p+\sum\limits_{j=2}^{k} (|\Delta a_j-\Delta c_j|^p+|\Delta b_j-\Delta d_j|^p) )^{\frac{1}{p}}$\\
		
		\hline

		Ashraf et al. \cite{2} & $S_{\rm Az_{p}}^{h}= 1-\frac{1}{\sqrt[p]{2k}}(|a_z-c_z|^p+|b_z-d_z|^p+\sum\limits_{j \in \delta(j)-z} (|\Delta a_j-\Delta c_j|^p+|\Delta b_j-\Delta d_j|^p)$\\&$~~~~~~~~~~+\sum\limits_{j \in \gamma(j)} (| a_j- c_j|^p+| b_j- d_j|^p) )^{\frac{1}{p}}$, where $\gamma(j) \cup \delta(j)= N,$ and $z=\mbox{min}\{\delta(j)\}$,\\&~~~~~~~~~~ $N$ is an index set of AIFS and $\gamma(j)~ \mbox{and}~ \delta(j)$ are two disjoint subset of $N$.\\ 
		
		
		\hline

		
		
		
		

	\end{tabular}
\end{table}

\section{Notations and Basic Definitions}\label{sec2}
Table \ref{nota} contains the list of symbols that are used in this paper. 
\begin{table}[htp]
	\centering
	\fontsize{8}{8}\selectfont
	\renewcommand{\arraystretch}{1.5}
	\caption{Notations}
	\label{nota}
	\begin{tabular}{p{2.0cm}|p{12.0cm}}
		\hline
		\textbf{Notations} & \textbf{Implications}\\
		\hline
		$y$ &  Sequence i.e. $(y_k)_{k=1}^\infty$, $y_k \in \mathbb{R}$ and $k \in \mathbb{N}$.\\ 
		\hline
		$x$ & Sequence i.e. $(x_k)_{k=1}^\infty$, $x_k \in \mathbb{R}$ and $k \in \mathbb{N}$.\\ 
		\hline
		$y-x$ &$(y_k - x_k)_{k=1}^{\infty}$, $k \in \mathbb{N}.$ \\
		\hline
		$\Delta y_k$&$(y_{k} - y_{k-1})_{k=1}^{\infty}, k \in \mathbb{N}.$ \\
		\hline
		$\Delta y$ &$(\Delta y_k)_{k=1}^{\infty}, k \in \mathbb{N}$ and provided $y_{0}=0.$  \\
		\hline
		$\tilde{y}$& Double Sequence i.e. $(y_{mn})_{m,n=1}^{\infty}$, $y_{mn} \in \mathbb{R}$ and $m,n \in \mathbb{N}.$\\
		\hline
		$\tilde{x}$& Double Sequence i.e. $(x_{mn})_{m,n=1}^{\infty}$, $x_{mn} \in \mathbb{R}$ and $m,n \in \mathbb{N}.$\\
		\hline
		$\tilde{y}-\tilde{x}$& $(y_{mn}-x_{mn})_{m,n=1}^{\infty,\infty}$ and $y_{mn}, x_{mn} \in \mathbb{R},~~ m,n \in \mathbb{N}.$\\
		\hline
		$\tilde{\underline{y}}_{}$& Two-dimensional AIFS sequence i.e., $(\mu_{k},\nu_{k})_{k=1}^{\infty}.$\\
		\hline
		$\tilde{\underline{x}}_{}$& Two-dimensional AIFS sequence i.e., $(\mu_{k}^{1},\nu_{k}^{1})_{k=1}^{\infty}.$\\
		\hline
		$ \tilde{\underline{y}}-\tilde{\underline{x}}$ & $\big(\mu_{k},\nu_{k})- (\mu_{k}^{1},\nu_{k}^{1})\big)_{k=1}^{\infty}.$	\\
		\hline
		$\Delta_{11} x_{mn}$& $x_{mn}-x_{mn-1}-x_{m-1n}+x_{m-1n-1}$ and $x_{00}=x_{10}=x_{01}=0.$\\
		\hline
		$\Delta_{11} \tilde{x}$& $(\Delta_{11}x_{mn})_{m,n=1}^{\infty,\infty}$ and $x_{00}=x_{10}=x_{01}=0.$\\
		\hline
		$\Omega$ & Collection of real-valued double sequences.\\
		\hline
		$\underline{\Omega}$ & Collection of two-dimensional real-valued sequences.\\
		\hline
		$\mathbb{R}$ & Set of real numbers.\\
		\hline
		$\mathbb{N}$ & Set of natural numbers.\\
		\hline
	\end{tabular}
	\label{table:summary MCDM}
\end{table}

\begin{definition}{Real valued double sequences} \cite{1}\\
	A function $\tilde{y}:\mathbb{N} \times \mathbb{N} \rightarrow \Omega$ defines a real valued double sequence such that $\tilde{y} = \{y_{mn}\}_{m,n=1}^{\infty}$, where 
	$y_{mn} \in \mathbb{R}, m,n \in \mathbb{N}.$ Here, $\Omega$ is the collection of real valued double sequences.
\end{definition}

\begin{definition}{Partition} \cite{2}\\
	In a real valued closed interval $[c,d]$, $k+1$ points $(\{y_0,y_1,y_2, \cdots y_k\})$ are selected for which $c=y_0 <y_1<y_2 \cdots<y_k=d$ holds, then these points results to a partition of $[c,d]$.
\end{definition}

\begin{definition}{Absolutely summable total variation} \cite{6}\\
	The absolutely summable total variation $(V_{\mbox{abs}})$ of a function $g$ over the partition $k$ is defined as:
	\begin{eqnarray}
		V_{\mbox{abs}}(g)=\bigg( |g(y_1)|+ \sum\limits_{i=2}^{k} |\Delta g(y_i)|\bigg)
	\end{eqnarray}
here, $\Delta g(y_i)= g(y_i)-g(y_{i-1}),$ for $i \geq 2.$
\end{definition}

\begin{definition}{Absolutely summable bounded variation} \cite{6}\\
	The absolutely summable bounded variation (BV$_{\mbox{abs}}$) is a collection of functions $g$ with a finite total variation. Mathematically,
		$	g \in BV_{\mbox{abs}([c, d])} ~\mbox{iff}~ \sup_{k} V_{\mbox{abs}}(g)	< \infty.$
\end{definition}

\begin{definition}{Real valued double sequences based metric space of absolutely summable bounded variation } \cite{1}\\
	A double sequence $\tilde{y}=\{y_{mn}\}_{m,n=1}^{\infty}$ which satisfies the Eq \ref{ab},
	\begin{eqnarray}\label{ab}
		\sum\limits_{m,n=1}^{\infty,\infty} |\Delta_{11}y_{mn}|^{}= \sum\limits_{m=1}^{\infty}\sum\limits_{n=1}^{\infty}	|\Delta_{11}y_{mn}|^{} < \infty
	\end{eqnarray}
	is termed as real valued double sequences of absolutely summable bounded variation, where
	$\Delta_{11} y_{mn}= y_{mn}-y_{m-1n}-y_{mn-1}+y_{m-1n-1}$. Here, $y_{m-1n}, y_{mn-1},$ and $y_{m-1n-1}$ is the three point neighbourhood of $y_{mn}$. In this space $	d(\tilde{y},\tilde{x})$ defines the metric as follows:
	\begin{eqnarray}
		d(\tilde{y},\tilde{x})=	\sum\limits_{m,n=1}^{\infty,\infty} |\Delta_{11}(y_{mn}-x_{mn})|^{}, 
	\end{eqnarray}
	$\mbox{~where,} ~x_{0n}=y_{m0}=x_{m0}=y_{0n}=y_{00}=x_{00}=0,~ m,n \in \mathbb{N}.$ 
\end{definition}



\begin{definition}{\normalfont \cite{4}}\normalfont Atanassov intuitionistic fuzzy set (AIFS) $B$ with respect to universe of discourse $U$ is defined as follows:
	\begin{equation}
			B=\left \{ (y,\mu_{B}(y),\nu_{B}(y)); y \in U \right \}
	\end{equation}
	Here, the membership function $\mu_{B}:U\rightarrow [0,1]$ and non-membership function $\nu_{B}:U\rightarrow [0,1]$ satisfies $0\leq \mu _{B}+\nu_{B}\leq 1$. In AIFS, $\pi_{B}=1-(\mu_{B}+\nu_{B})$ is the hesitancy function. 
\end{definition}
\begin{Remark}
	We have considered the two components $(\mu_{B}(y),\nu_{B}(y))=(\mu,\nu) $ of AIFS as independent variables. The two-dimensional plane is used to depict $(\mu,\nu)$. This plane is termed an AIFS plane because it captures the two-dimensional view of AIFS.
\end{Remark}

\begin{definition}{Two-dimensional AIFS sequences}\\
	A function $\tilde{\underline{y}}: \rm AIFS(U) \rightarrow \underline{\Omega}$ defines a two-dimensional AIFS sequence such that $\tilde{\underline{y}}=(\mu_{k},\nu_{k})_{k=1}^{\infty}$.
Here, $\underline{\Omega}$ is the collection of two-dimensional real valued sequences.
\end{definition}

\begin{definition}
	\cite{14}
	Let $B_1$, $B_2$, and $B_3$ be AIFSs belonging to an AIFS plane AIFS(U). A measure $D: \rm AIFS(U) \times \rm AIFS(U)$ $\rightarrow$ [0,1] satisfying the following four properties is known as the distance measure. 
	\begin{enumerate}    
		\item $0 \leq D(B_1,B_2) \leq 1$
		\item $D(B_1,B_2)$ = 0 $\iff$ $B_1$ = $B_2$
		\item $D(B_1,B_2)$ = $D(B_2,B_1)$    
		\item If $B_1 \subseteq B_2 \subseteq B_3$, then $D(B_1,B_3)$ $\geq$ $D(B_1,B_2)$ and $D(B_1,B_3)$  $\geq$ $D(B_2,B_3)$
	\end{enumerate}
	The set $(U, D)$ together is known as metric space.
\end{definition}

\begin{definition}
	\cite{19}
	The similarity measure is a mapping $S: \rm AIFS(U) \times \rm AIFS(U) \rightarrow$ [0,1], which satisfies the underneath four properties. Let $B_1$, $B_2$, and $B_3$ be the AIFSs belonging to AIFS(U).
	\begin{enumerate}
		\item $0 \leq S(B_1,B_2) \leq 1$
		\item  $ S(B_1,B_2)$ = 1 $\iff B_1 = B_2$
		\item $S(B_1,B_2) = S(B_2,B_1)$
		\item If $B_1 \subseteq B_2 \subseteq B_3$, then $S(B_1,B_3) \leq S(B_1,B_2)$ and $S(B_1,B_3) \leq S(B_2,B_3)$
	\end{enumerate}
\end{definition}

Similarity and distance measure are the complement of each other, and this complementary relationship is  
$S(B_1,B_2)=1-D(B_1,B_2).$ 



\section{Main Results}\label{sec3}
This section has two subsections. In subsection-\ref{3.1}, some definitions are proposed for the introduction of bounded variation based spatial metric space. The spatial metric space, its completeness, and spatial similarity measure are proposed in subsection-\ref{3.2}. 
\subsection{Novel divided difference operator}\label{3.1}
In this section, an overview of the difference operator and novel divided difference operator is given in terms of the geometrical interpretation. Further, the distributive property of these operators are also studied besides defining a novel divided difference operator.
\begin{definition}
	{Geometrical interpretation of difference operator $\Delta$}\\
	Let $y=(y_{k})_{k=1}^{\infty}$ be a real-valued sequence. The step sizes of $y$ are estimated by the help of ordinary difference operator $\Delta$, provided $y_{0}=0$ is equal to
	$$\Delta y=(\Delta y_k)_{k=1}^{\infty}=(y_{k} - y_{k-1})_{k=1}^{\infty}$$
	The difference of the point $y_{k}$ with the point $y_{k-1}$ gives the step size involved in moving from $y_{k-1}$ to its succeeding point $y_{k}$. The value of $\Delta y$ informs about the variation in the step size. 
	In sequence succeeding element act as a neighbourhood of the preceding element;
	hence every element has one point based neighbourhood.
	We have, $$\Delta(y-x)=\Delta y-\Delta x$$
\end{definition} 
Now, we study the distributive property of difference operator $\Delta$ over a sequence of two-dimensional AIFS,  $\tilde{\underline{y}}=(\mu_{k},\nu_{k})_{k=1}^{\infty}$ as follows:
\begin{Theorem}
	Let $\tilde{\underline{y}}$ and $\tilde{\underline{x}}$ are the sequences of two dimensional AIFS, we have
	$\Delta(\tilde{\underline{y}}-\tilde{\underline{x}})=\Delta \tilde{\underline{y}}-\Delta \tilde{\underline{x}}$.
\end{Theorem}
\begin{Proof}
	Let $\tilde{\underline{y}}=(\mu_{k},\nu_{k})_{k=1}^{\infty}$ and $\tilde{\underline{x}}=(\mu^{1}_{k},\nu^{1}_{k})_{k=1}^{\infty}$ such that
	and $\mu_{k},\mu^{1}_{k},\nu_{k},\nu^{1}_{k}
	\in [0,1]$. Now, $$\tilde{\underline{y}}-\tilde{\underline{x}}=(\mu_{k}-\mu^{1}_{k},\nu_{k}-\nu^{1}_{k})_{k=1}^{\infty}$$. Further, $$\Delta \tilde{\underline{y}}= (\Delta \mu_k, \Delta \nu_k)_{k=1}^{\infty}=(\mu_{k}-\mu_{k-1},\nu_{k}-\nu_{k-1})_{k=1}^{\infty}$$
	We have,
	\begin{eqnarray}
		\Delta (\tilde{\underline{y}}-\tilde{\underline{x}})&=& ((\mu_{k}-\mu^{1}_{k})- (\mu_{k-1}-\mu^{1}_{k-1}), (\nu_{k}-\nu^{1}_{k})- (\nu_{k-1}-\nu^{1}_{k-1}))_{k=1}^{\infty} \nonumber\\
		&=& ((\mu_{k}-\mu_{k-1})- (\mu^{1}_{k}-\mu^{1}_{k-1}), (\nu_{k}-\nu_{k-1})- (\nu^{1}_{k}-\nu^{1}_{k-1}))_{k=1}^{\infty} \nonumber\\
		&=&\Delta \tilde{\underline{y}}  - \Delta \tilde{\underline{x}}. \nonumber ~~~~~~~~~~~~~\square
	\end{eqnarray}
\end{Proof}
\begin{definition}{Divided difference operator $\Delta_{11}$} \cite{25}\\
	The divided difference operator $\Delta_{11}$  is defined on a real-valued double sequences $\widetilde {y}=(y_{mn})_{m,n=1}^{\infty}$ as follows:
	$$\Delta_{11}y_{mn}=y_{mn}-y_{mn-1}-y_{m-1n}+y_{m-1n-1}$$
	The step sizes in double sequences $\widetilde {y}$ are estimated with the help of divided difference operator $\Delta_{11}$.
\end{definition}
\begin{definition}{Novel divided difference operator $\Delta_{11}^{1}$}\\
	Let $\tilde{\underline{y}}=(\mu_{k},\nu_{k})_{k=1}^{\infty}$ be a sequence of two-dimensional AIFS. In $\tilde{\underline{y}}; ~\nu_k, ~\mu_{k-1},$ and $\nu_{k-1}$ are
	clubbed to $\mu_{k}$ in $(\mu_k,\nu_k), (\mu_{k-1},\nu_{k-1})$. So $\nu_k, ~\mu_{k-1},$ and $\nu_{k-1}$ is a three-point based neighbourhood of $\mu_{k}$ with respect to $(\mu_k,\nu_k)$, $(\mu_{k-1},\nu_{k-1})$. Now, $\Delta_{11}^{1}$ is defined as follows: $$\Delta_{11}^{1}((\mu_{k},\nu_{k}),(\mu_{k-1},\nu_{k-1}))=\mu_{k}-\mu_{k-1}+\nu_{k-1}-\nu_{k}$$ provided that $\mu_{0}=\nu_{0}=0.$
\end{definition}


\begin{definition}{Geometrical interpretation of novel divided difference operator $\Delta_{11}^{1}$}\\
	The novel divided difference operator operates on $(\mu_k,\nu_k)$, $(\mu_{k-1},\nu_{k-1})$ as they possess a three point based neighbourhood. Here, a unit length is decreased  along the membership and non-membership components while moving a step. So, the movement from $(\mu_{k-1},\nu_{k-1})$ to $(\mu_k,\nu_k)$ in AIFS plane involves a move from $\mu_{k}$ to $\mu_{k-1}$ and another move from $\nu_{k-1}$ to $\nu_{k}$. The difference between the
	movements from $\mu_{k}$ to $\mu_{k-1}$ and $\nu_{k}$ to $\nu_{k-1}$ gives the step size, this difference is calculated with the help of $\Delta_{11}^{1}((\mu_k,\nu_k)$, $(\mu_{k-1},\nu_{k-1}))$. The variation in the step sizes is shown in the components of $\Delta_{11}^{1} \tilde{{\underline{y}}}.$
\end{definition}
\begin{Theorem}\label{thm3.2}
	Let $\tilde{\underline{y}}$ and $\tilde{\underline{x}}$ are the sequences of two dimensional AIFS, we have $\Delta_{11}^{1}(\tilde{\underline{y}}-\tilde{\underline{x}})=\Delta_{11}^{1}\tilde{\underline{y}}-\Delta_{11}^{1}\tilde{\underline{x}}.$
\end{Theorem}
\begin{Proof}
	Let $\tilde{\underline{y}}=(\mu_{k},\nu_{k})_{k=1}^{\infty}$ and $\tilde{\underline{x}}=(\mu^{1}_{k},\nu^{1}_{k})_{k=1}^{\infty}$ such that
	and $\mu_{k},\mu^{1}_{k},\nu_{k}, \nu^{1}_{k}
	\in [0,1]$. Now, $$\tilde{\underline{y}}-\tilde{\underline{x}}=(\mu_{k}-\mu^{1}_{k},\nu_{k}-\nu^{1}_{k})_{k=1}^{\infty}$$. Further, $$\Delta_{11}^{1}\underline{{y}}= \big(\Delta_{11}^{1}((\mu_k,\nu_k), (\mu_{k-1},\nu_{k-1})) \big)_{k=1}^{\infty}=(\mu_{k}-\mu_{k-1}+\nu_{k-1}-\nu_{k})_{k=1}^{\infty}$$
	We have,
	\begin{eqnarray}
		\Delta_{11}^{1} (\tilde{\underline{y}}-\tilde{\underline{x}})&=& (\mu_{k}-\mu^{1}_{k})- (\mu_{k-1}-\mu^{1}_{k-1})+(\nu_{k-1}-\nu^{1}_{k-1})-(\nu_{k}-\nu^{1}_{k})\nonumber\\
		&=& (\mu_{k}-\mu_{k-1}+\nu_{k-1}-\nu_{k})-(\mu^{1}_{k}-\mu^{1}_{k-1}+\nu^{1}_{k-1}-\nu^{1}_{k}) \nonumber\\
		&=&\Delta_{11}^{1} \tilde{\underline{y}}  - \Delta_{11}^{1} \tilde{\underline{x}}.  \nonumber ~~~~~~~~~~~~~\square
	\end{eqnarray}
\end{Proof}
\subsection{Spatial metric space of bounded variation and spatial similarity measure} \label{3.2}
This section defines a novel divided difference operator based on an absolutely summable bounded variation of two-dimensional AIFS valued sequences. A mathematical study of spatial metric space and its convergence is carried out to establish the proposed spatial similarity measure. 
\begin{definition}{Multivalued absolutely summable bounded variation in AIFS(U)}\\
	The absolutely summable bounded variation in AIFS(U) is a two-dimensional AIFS valued sequence $\tilde{\underline{y}}=(\mu_k,\nu_k)_{k=1}^{\infty}$ which satisfies the following conditions:
	\begin{eqnarray}
	&1.&	\sum_{k=1}^{\infty}	|\Delta_{11}^{1}(({\mu}_{k},0),({\mu}_{k-1},0))|=\sum_{k=1}^{\infty} |\mu_{k}-\mu_{k-1}| < \infty, \nonumber \\
	&2.&		\sum_{k=1}^{\infty}	|\Delta_{11}^{1}(({\mu}_{k},\nu_{k}),({\mu}_{k-1},\nu_{k-1}))|=\sum_{k=1}^{\infty} |\mu_{k}-\mu_{k-1}+\nu_{k-1}-\nu_{k}| < \infty \nonumber
	\end{eqnarray}
where, $\mu_{0}=\nu_{0}=0.$ Here, $\mu_{k-1}$, $\nu_{k-1}$, and $\nu_{k}$ is the three point neighbourhood of $\mu_{k}$ with respect to $(\mu_k,\nu_k),(\mu_{k-1},\nu_{k-1})$.
\end{definition}

\begin{Theorem} \label{3.3}
    Let $\underline{y}^{k}$ be a collection of sequence of two-dimensional AIFS. For $\tilde{{\underline{y}}} \in \underline{y}^{k}$, we have $\tilde{{\underline{y}}}=(\mu_{j},\nu_{j})_{j=1}^{\infty}$ such that $\mu_{j},\nu_{j}\geq 0,~\forall~ j=1,2,\cdots k,$ and $\mu_{j+1},\nu_{j+1}=0,~\forall~j>k.$ The distance measure over $\underline{y}^{k}$  simulatenously involves three distance measure $D_1,D_2$, and $D_3$. Here,
     $D_1: \underline{y}^{k} \times \underline{y}^{k} \rightarrow [0,1]$, $D_2: \underline{y}^{k} \times \underline{y}^{k} \rightarrow [0,1]$ and $D_3: \underline{y}^{k} \times \underline{y}^{k} \rightarrow [0,1]$ such that
\begin{eqnarray}
	D_1(\tilde{{\underline{y}}},\tilde{{\underline{x}}})&=&	 \frac{1}{4k}\sum\limits_{j=1}^{k}\bigg(\big|\Delta_{11}^{1}\big((\mu_j,0),(\mu_{j-1},0)\big)
	-\Delta_{11}^{1}\big((\mu_j^{1},0),(\mu_{j-1}^{1},0)\big)\big|  \nonumber \\
&&\hspace{2.0cm}+\big|\Delta_{11}^{1}\big((\mu_j,\nu_{j}),(\mu_{j-1},\nu_{j-1})\big)
-\Delta_{11}^{1}\big((\mu_j^{1},\nu_j^{1}),(\mu_{j-1}^{1},\nu_{j-1}^{1})\big)\big|\bigg) \label{41} \\
D_2(\tilde{{\underline{y}}},\tilde{{\underline{x}}})&=&	 \frac{1}{4k}\sum\limits_{j=1}^{k}\bigg(\big|\Delta_{11}^{1}\big((0,\nu_j),(0,\nu_{j-1})
\big)
-\Delta_{11}^{1}\big((0,\nu_j^{1}),(0,\nu_{j-1}^{1})\big)\big| \nonumber\\
&&\hspace{2.0cm}+\big|\Delta_{11}^{1}\big((\mu_j,\nu_j),(\mu_{j-1},\nu_{j-1})\big)
-\Delta_{11}^{1}\big((\mu_j^{1},\nu_j^{1}),(\mu_{j-1}^{1},\nu_{j-1}^{1})\big)\big|\bigg) \label{51}	 \\
D_3(\tilde{{\underline{y}}},\tilde{{\underline{x}}})&=&	 \frac{1}{4k}\sum\limits_{j=1}^{k}\bigg(\big|\Delta_{11}^{1}\big((\mu_j,0),(\mu_{j-1},0)\big)
-\Delta_{11}^{1}\big((\mu_j^{1},0),(\mu_{j-1}^{1},0)\big) \big|  \nonumber\\
&&\hspace{2.0cm}+\big|\Delta_{11}^{1}\big((0,\nu_j),(0,\nu_{j-1})
\big)-\Delta_{11}^{1}\big((0,\nu_j^{1}),(0,\nu_{j-1}^{1})\big)\big| \bigg) \label{6}
\end{eqnarray}
We have termed $D_1$ as the membership dominant distance measure, $D_2$ as the non-membership dominant distance measure, and $D_3$ as the equidominant distance measure over $\underline{y}^{k}$. It is given that $\mu_{0}=\nu_{0}=\mu_{0}^{1}= \nu_{0}^{1}=0$. The collection $\bigg(\underline{y}^{k},(D_1,D_2,D_3)\bigg)$ is the spatial metric space.
\end{Theorem}

\begin{Proof}
Let us verify the four properties of a distance measure for $D_1$.
\begin{itemize}
	\item[1.] To show, $0 \leq 	D_1(\tilde{{\underline{y}}},\tilde{{\underline{x}}}) \leq 1$. 
	\begin{eqnarray}	
		&& \mbox{Since,~~~}0 \leq \bigg(\sum\limits_{j=1}^{k}\big|\Delta_{11}^{1}\big((\mu_j,\nu_j),(\mu_{j-1},\nu_{j-1})\big)-\Delta_{11}^{1}\big((\mu_j^{1},\nu_j^{1}),(\mu_{j-1}^{1},\nu_{j-1}^{1})\big)\big|\bigg)  \leq 4(k-1)+2 \label{4} \\
		&& \mbox{and,} ~~~0 \leq \bigg(\sum\limits_{j=1}^{k}\big|\Delta_{11}^{1}\big((\mu_j,0),(\mu_{j-1},0)\big)-\Delta_{11}^{1}\big((\mu_j^{1},0),(\mu_{j-1}^{1},0)\big)\big|\bigg) \leq 2 \label{5} 
	\end{eqnarray} 
	Adding, Eq. \ref{4} and \ref{5}, we get
	\begin{eqnarray}
		&&  0 \leq \bigg(\sum\limits_{j=1}^{k}\big|\Delta_{11}^{1}\big((\mu_j,0),(\mu_{j-1},0)\big)-\Delta_{11}^{1}\big((\mu_j^{1},0),(\mu_{j-1}^{1},0)\big)\big|\bigg) \nonumber\\
		&&\hspace*{1.7cm} +\bigg(\sum\limits_{j=1}^{k}\big|\Delta_{11}^{1}\big((\mu_j,\nu_j),(\mu_{j-1},\nu_{j-1})\big)-\Delta_{11}^{1}\big((\mu_j^{1},\nu_j^{1}),(\mu_{j-1}^{1},\nu_{j-1}^{1})\big)\big|\bigg) \leq 4(k-1)+2+2\leq 4k \nonumber\\
		&&\hspace*{0.0cm}\mbox{Hence,~}~~0 \leq  {\frac{1}{4k} \sum\limits_{j=1}^{k}\bigg(\big|\Delta_{11}^{1}\big((\mu_j,0),(\mu_{j-1},0)\big)-\Delta_{11}^{1}\big((\mu_j^{1},0),(\mu_{j-1}^{1},0)\big)\big|     } \nonumber\\
		&&\hspace*{3.5cm}+\big|\Delta_{11}^{1}\big((\mu_j,\nu_i),(\mu_{j-1},\nu_{j-1})\big)-\Delta_{11}^{1}\big((\mu_j^{1},\nu_j^{1}),(\mu_{j-1}^{1},\nu_{j-1}^{1})\big)\big|\bigg) \leq 1 \nonumber
	\end{eqnarray}
	$\implies 0 \leq 	D_1(\tilde{{\underline{y}}},\tilde{{\underline{x}}}) \leq 1.$
	
	\item[2.] To show,  ${\underline{y}}={\underline{x}} \mbox{~if and only if~} ~D_1(\tilde{{\underline{y}}},\tilde{{\underline{x}}})$ = 0.\\
	Let $	D_1(\tilde{{\underline{y}}},\tilde{{\underline{x}}})=0$,
	\begin{eqnarray}
		&&\Longleftrightarrow {\frac{1}{4k} \sum\limits_{j=1}^{k}\bigg(\big|\Delta_{11}^{1}\big((\mu_j,0),(\mu_{j-1},0)\big)-\Delta_{11}^{1}\big((\mu_j^{1},0),(\mu_{j-1}^{1},0)\big)\big|     } \nonumber\\
		&&\hspace*{2.5cm}+\big|\Delta_{11}^{1}\big((\mu_j,\nu_j),(\mu_{j-1},\nu_{j-1})\big)-\Delta_{11}^{1}\big((\mu_j^{1},\nu_j^{1}),(\mu_{j-1}^{1},\nu_{j-1}^{1})\big)\big|\bigg)  =0. \nonumber
	\end{eqnarray}
	$\hspace*{2.5cm}\Longleftrightarrow \big|\Delta_{11}^{1}\big((\mu_j,0),(\mu_{j-1},0)\big)-\Delta_{11}^{1}\big((\mu_j^{1},0),(\mu_{j-1}^{1},0)\big)\big|=0, $ and\\
	$\hspace*{3.5cm}\big|\Delta_{11}^{1}\big((\mu_j,\nu_j),(\mu_{j-1},\nu_{j-1})\big)-\Delta_{11}^{1}\big((\mu_j^{1},\nu_j^{1}),(\mu_{j-1}^{1},\nu_{j-1}^{1})\big)\big|=0.$
	
	Now, $\big|\Delta_{11}^{1}\big((\mu_j,0),(\mu_{j-1},0)\big)-\Delta_{11}^{1}\big((\mu_j^{1},0),(\mu_{j-1}^{1},0)\big)\big|=0$, implies
	\begin{eqnarray}\label{8}
		&&|(\mu_1-\mu_0)-(\mu_1^{1}-\mu_0^{1})|=0 \implies \mu_1=\mu_1^{1}  \nonumber\\
		&& |(\mu_2-\mu_1)-(\mu_2^{1}-\mu_1^{1})|=0 \implies \mu_2=\mu_2^{1}  \nonumber\\
		&& \vdots \nonumber \\
		&& |(\mu_k-\mu_{k-1})-(\mu_k^{1}-\mu_{k-1}^{1})|=0 \implies \mu_k=\mu_k^{1}  
	\end{eqnarray}
	From $\big|\Delta_{11}^{1}\big((\mu_j,\nu_j),(\mu_{j-1},\nu_{j-1})\big)-\Delta_{11}^{1}\big((\mu_j^{1},\nu_j^{1}),(\mu_{j-1}^{1},\nu_{j-1}^{1})\big)\big|=0$, we get 
	\begin{eqnarray}\label{9}
		&&|(\mu_1-\mu_0+\nu_0-\nu_1)-(\mu_1^{1}-\mu_0^{1}+\nu_0^{1}-\nu_1^{1})|=0 \implies \nu_1=\nu_1^{1}  \nonumber\\
		&&|(\mu_2-\mu_1+\nu_1-\nu_2)-(\mu_2^{1}-\mu_1^{1}+\nu_1^{1}-\nu_2^{1})|=0 \implies \nu_2=\nu_2^{1}  \nonumber\\
		&& \vdots \nonumber\\ &&\hspace*{-1.5cm}|(\mu_k-\mu_{k-1}+\nu_{k-1}-\nu_k)-(\mu_k^{1}-\mu_{k-1}^{1}+\nu_{k-1}^{1}-\nu_{k}^{1})|=0 \implies \nu_k=\nu_k^{1} \nonumber 
	\end{eqnarray}
	$\Longleftrightarrow \tilde{{\underline{y}}}=\tilde{{\underline{x}}}$.
	
	\item[3.] $	D_1(\tilde{{\underline{y}}},\tilde{{\underline{x}}})=	D_1(\tilde{{\underline{x}}},\tilde{{\underline{y}}})$ is obvious.
	\item[4.] To show, If $\tilde{{\underline{y}}} \subseteq \tilde{{\underline{x}}} \subseteq \tilde{{\underline{z}}}$, where $\tilde{{\underline{y}}}$,$\tilde{{\underline{x}}}$,$\tilde{{\underline{z}}}$ $\in$ $AIFSs(U)$, then $D_1(\tilde{{\underline{y}}},\tilde{{\underline{z}}})$ $\geq$ $D_1(\tilde{{\underline{y}}},\tilde{{\underline{x}}})$ and $D_1(\tilde{{\underline{y}}},\tilde{{\underline{z}}})$  $\geq$ $D_1(\tilde{{\underline{x}}},\tilde{{\underline{z}}})$.\\
	Instead of this property, we will prove the triangle inequality. Since,
	\begin{eqnarray}
			D_1(\tilde{{\underline{y}}},\tilde{{\underline{x}}})&=& {\frac{1}{4k} \sum\limits_{j=1}^{k}\bigg(\big|\Delta_{11}^{1}\big((\mu_j,0),(\mu_{j-1},0)\big)-\Delta_{11}^{1}\big((\mu_j^{1},0),(\mu_{j-1}^{1},0)\big)\big|     } \nonumber\\
		&&\hspace*{2.5cm}+\big|\Delta_{11}^{1}\big((\mu_j,\nu_j),(\mu_{j-1},\nu_{j-1})\big)-\Delta_{11}^{1}\big((\mu_j^{1},\nu_j^{1}),(\mu_{j-1}^{1},\nu_{j-1}^{1})\big)\big|\bigg) \nonumber\\
		&=& \frac{1}{4k} \hspace*{0.0cm}\sum\limits_{j=1}^{k}\bigg(\big|\Delta_{11}^{1}\big((\mu_j,0),(\mu_{j-1},0)\big)-\Delta_{11}^{1}\big((\mu_j^{2},0),(\mu_{j-1}^{2},0)\big)\nonumber\\
		&&\hspace*{2.5cm} +\Delta_{11}^{1}\big((\mu_j^{2},0),(\mu_{j-1}^{2},0)\big)-\Delta_{11}^{1}\big((\mu_j^{1},0),(\mu_{j-1}^{1},0)\big)\big|\nonumber\\
		&&	\hspace*{2.5cm}+\big|\Delta_{11}^{1}\big((\mu_j,\nu_j),(\mu_{j-1},\nu_{j-1})\big)
		-\Delta_{11}^{1}\big((\mu_j^{2},\nu_j^{2}),(\mu_{j-1}^{2},\nu_{j-1}^{2})\big)\nonumber\\
		&&\hspace*{2.5cm}+ \Delta_{11}^{1}\big((\mu_j^{2},\nu_j^{2}),(\mu_{j-1}^{2},\nu_{j-1}^{2})\big)-\Delta_{11}^{1}\big((\mu_j^{1},\nu_j^{1}),(\mu_{j-1}^{1},\nu_{j-1}^{1})\big)\big|\bigg) \nonumber\\
		&\leq & \frac{1}{4k} \hspace*{0.0cm}\sum\limits_{j=1}^{k}\bigg(\big|\Delta_{11}^{1}\big((\mu_j,0),(\mu_{j-1},0)\big)-\Delta_{11}^{1}\big((\mu_j^{2},0),(\mu_{j-1}^{2},0)\big)\big|\nonumber\\
		&& \hspace*{2.5cm}+\big|\Delta_{11}^{1}\big((\mu_j^{2},0),(\mu_{j-1}^{2},0)\big)-\Delta_{11}^{1}\big((\mu_j^{1},0),(\mu_{j-1}^{1},0)\big)\big|\nonumber\\
		&&	\hspace*{2.5cm}+\big|\Delta_{11}^{1}\big((\mu_j,\nu_j),(\mu_{j-1},\nu_{j-1})\big)
		-\Delta_{11}^{1}\big((\mu_j^{2},\nu_j^{2}),(\mu_{j-1}^{2},\nu_{j-1}^{2})\big)\nonumber\\
	&&\hspace*{2.5cm}+\Delta_{11}^{1}\big((\mu_j^{2},\nu_j^{2}),(\mu_{j-1}^{2},\nu_{j-1}^{2})\big)-\Delta_{11}^{1}\big((\mu_j^{1},\nu_j^{1}),(\mu_{j-1}^{1},\nu_{j-1}^{1})\big)\big|\bigg) \nonumber\\
		\hspace*{1.0cm}&=& \frac{1}{4k} \sum\limits_{j=1}^{k}\bigg(\big|\Delta_{11}^{1}\big((\mu_j,0),(\mu_{j-1},0)\big)-\Delta_{11}^{1}\big((\mu_j^{2},0),(\mu_{j-1}^{2},0)\big)\big|\nonumber\\
		&&\hspace*{2.5cm}+
		\Delta_{11}^{1}\big((\mu_j,\nu_j),(\mu_{j-1},\nu_{j-1})\big)-\Delta_{11}^{1}\big((\mu_j^{2},\nu_j^{2}),(\mu_{j-1}^{2},\nu_{j-1}^{2})\big)\big|\bigg) \nonumber\\
	&& \hspace*{2.5cm}+\frac{1}{4k} \hspace*{0.0cm}\sum\limits_{j=1}^{k}\bigg(\big|\Delta_{11}^{1}\big((\mu_j^{2},0),(\mu_{j-1}^{2},0)\big)-\Delta_{11}^{1}\big((\mu_j^{1},0),(\mu_{j-1}^{1},0)\big)\big|\nonumber\\
	&&\hspace*{2.5cm}+\Delta_{11}^{1}\big((\mu_j^{2},\nu_j^{2}),(\mu_{j-1}^{2},\nu_{j-1}^{2})\big)-\Delta_{11}^{1}\big((\mu_j^{1},\nu_j^{1}),(\mu_{j-1}^{1},\nu_{j-1}^{1})\big)\big|\bigg) \nonumber\\
		&=& 	D_1(\tilde{{\underline{y}}},\tilde{{\underline{z}}})+	D_1(\tilde{{\underline{z}}},\tilde{{\underline{x}}}) \nonumber 
	\end{eqnarray}
	In the similar fashion $D_2$ and $D_3$ can be proved to be a distance measure.\\
	Hence, there exists three distance measures ($D_1,D_2,D_3$) over $\underline{y}^k$ simulatenously, and it results to a spatial metric space $\bigg(\underline{y}^{k},(D_1,D_2,D_3)\bigg)$. $~~~~~~~~~~~~~\square $ 
\end{itemize}
\end{Proof}
\begin{Remark}
	The factor $\big|\Delta_{11}^{1}\big((\mu_j,0),(\mu_{j-1},0)\big)
	-\Delta_{11}^{1}\big((\mu_j^{1},0),(\mu_{j-1}^{1},0)\big)\big|$ gives an additional importance to membership values in $D_1$ and hence $D_1$ is referred as membership dominant distance measure. Similarly, the factor $\big|\Delta_{11}^{1}\big((0,\nu_j),(0,\nu_{j-1})
	\big)
	-\Delta_{11}^{1}\big((0,\nu_j^{1}),(0,\nu_{j-1}^{1})\big)\big| $ gives an additional importance to non-membership values in $D_2$ and hence it is termed as non-membership dominant distance measure. In $D_3$, the factors $\big|\Delta_{11}^{1}\big((\mu_j,0),(\mu_{j-1},0)\big)
	-\Delta_{11}^{1}\big((\mu_j^{1},0),(\mu_{j-1}^{1},0)\big)\big|$ and $\big|\Delta_{11}^{1}\big((0,\nu_j),(0,\nu_{j-1})
	\big)
	-\Delta_{11}^{1}\big((0,\nu_j^{1}),(0,\nu_{j-1}^{1})\big)\big| $  are giving importance to membership and non-membership values individually, so it is referred as equidominant distance measure.
\end{Remark}
\textbf{Corollary 3.3.1.}
 Let $\underline{y'}^{k}$ be a collection of sequence of fuzzy sets upto $k$ indices that is, $\tilde{{\underline{y}}} \in \underline{y'}^{k}$, we have $\tilde{{\underline{y}}}=(\mu_{j},\nu_{j})_{j=1}^{\infty}$ such that $\mu_{j}$,$\nu_{j}\geq 0, ~\mu_{j}+\nu_{j}=1~\forall~ j=1,2,\cdots k,$ and $\mu_{j+1},\nu_{j+1}=0,~\forall~j>k.$
Here, $D_1,D_2$, and $D_3$ reduce to $D_1'$ and $D_2'$, where
\begin{eqnarray}
	D_1'(\tilde{{\underline{y}}},\tilde{{\underline{x}}})&=&	 \frac{1}{4k}\sum\limits_{j=1}^{k}\bigg(\big|\Delta_{11}^{1}\big((\mu_j,0),(\mu_{j-1},0)\big)
	-\Delta_{11}^{1}\big((\mu_j^{1},0),(\mu_{j-1}^{1},0)\big)\big|\bigg)  \\
	D_2'(\tilde{{\underline{y}}},\tilde{{\underline{x}}})&=&	 \frac{1}{2k}\sum\limits_{j=1}^{k}\bigg(\big|\Delta_{11}^{1}\big((\mu_j,0),(\mu_{j-1},0)\big)
	-\Delta_{11}^{1}\big((\mu_j^{1},0),(\mu_{j-1}^{1},0)\big)\big|\bigg)
\end{eqnarray}
In this case, $D_1$ and $D_2$ becomes $D_1'$, which results a two branched multivalued distance measure.\\
\textbf{Proof 3.3.1.} In a fuzzy set, we have $\mu_{j}=1-\nu_{j},$ and thus,
   \begin{eqnarray}
   &&\Delta_{11}^{1}\big((\mu_j,\nu_{j}),(\mu_{j-1},\nu_{j-1})\big)
    -\Delta_{11}^{1}\big((\mu_j^{1},\nu_j^{1}),(\mu_{j-1}^{1},\nu_{j-1}^{1})\big)=0. \label{12}\\
    &&\Delta_{11}^{1}\big((\mu_j,0),(\mu_{j-1},0)\big)
    -\Delta_{11}^{1}\big((\mu_j^{1},0),(\mu_{j-1}^{1},0)\big) =
    \Delta_{11}^{1}\big((0,\nu_j),(0,\nu_{j-1})
    \big)-\Delta_{11}^{1}\big((0,\nu_j^{1}),(0,\nu_{j-1}^{1})\big)\label{13}
\end{eqnarray}
Now, Eq.\ref{41},\ref{51},\ref{12} results $D_1'$ and from Eq.\ref{6},\ref{13} $D_2'$ is obtained. $\hspace*{2.0cm}\square$\\
~\\
\textbf{Corollary 3.3.2.} $\mathbb{R}$ is equivalent to $\mathbb{R} \times 0$ and $0 \times\mathbb{R}.$ 
Let $\underline{y''}^{k}$ be a collection of real-valued bounded variation sequences. For $\tilde{{\underline{y}}} \in \underline{y''}^{k}$, we have $\tilde{{\underline{y}}}=(y_{j})_{j=1}^{\infty}$
such that $0 \leq y_{j} \leq 1,~\forall~ j=1,2,\cdots k,$ and $y_{j+1}=0,~\forall~j>k.$ Moreover $y_{j}$ is equivalent to $(y_{j}, 0)$ or $(0, y_{j})$. Here, $D_1,D_2$, and $D_3$ reduces to a univalued distance measure $D,$ where, $D: \underline{y''}^{k} \times \underline{y''}^{k} \rightarrow [0,1]$ such that
\begin{eqnarray}
	D(\tilde{{\underline{y}}},\tilde{{\underline{x}}})&=&	 \frac{1}{2k}\sum\limits_{j=1}^{k}\bigg(\big|\Delta_{11}^{1}\big((y_j,0),(y_{j-1},0)\big)
	-\Delta_{11}^{1}\big((x_j^{1},0),(x_{j-1}^{1},0)\big)\big|\bigg) 
\end{eqnarray}
\hspace*{-0.0cm}\textbf{Proof 3.3.2.} Since, $\mathbb{R}$ is equivalent to $\mathbb{R} \times 0$ and $0 \times\mathbb{R}.$ Therefore, the operator $\Delta_{11}^{1}$ is calculateed over $\mathbb{R}$ in two ways as follows:  
\begin{eqnarray}
	\Delta_{11}^{1}\big((0,x_j),(0,x_{j-1})
	\big)-\Delta_{11}^{1}\big((0,x_j^{1}),(0,x_{j-1}^{1})\big)
\mbox{~~or,~~}	\Delta_{11}^{1}\big((y_j,0),(y_{j-1},0)\big)-\Delta_{11}^{1}\big((y_j^{1},0),(y_{j-1}^{1},0)\big) \nonumber
\end{eqnarray} 
for $y_j, x_j\in\mathbb{R}$. Now,\\
$\Delta_{11}^{1}\big((y_j,x_{j}),(y_{j-1},x_{j-1})\big)
-\Delta_{11}^{1}\big((y_j^{1},x_j^{1}),(y_{j-1}^{1},x_{j-1}^{1})\big)$
\begin{eqnarray}\label{15}
	&=& \left\{
	\begin{array}{lr}
			\Delta_{11}^{1}\big((y_j,0),(y_{j-1},0)\big)-\Delta_{11}^{1}\big((y_j^{1},0),(y_{j-1}^{1},0)\big), & \mbox{along~} \mathbb{R} \times 0  \\
		\Delta_{11}^{1}\big((0,x_j),(0,x_{j-1})
		\big)-\Delta_{11}^{1}\big((0,x_j^{1}),(0,x_{j-1}^{1})\big),  & \mbox{along}~0 \times\mathbb{R}  
	  \end{array}
	\right.
\end{eqnarray}
We have, \begin{eqnarray} \label{16}
	\Delta_{11}^{1}\big((y_j,0),(y_{j-1},0)\big)-\Delta_{11}^{1}\big((y_j^{1},0),(y_{j-1}^{1},0)\big)=\Delta_{11}^{1}\big((0,x_j),(0,x_{j-1})
	\big)-\Delta_{11}^{1}\big((0,x_j^{1}),(0,x_{j-1}^{1})\big)
\end{eqnarray}
using Eq.\ref{15},\ref{16} in Eq.\ref{41},\ref{51},\ref{6}, we get $D.$ $\hspace*{2.0cm}\square$

\begin{Remark}
	In the concept of bounded variation, a novel divided difference is proposed in the fuzzy/two-dimensional AIFS domain such that the notion of multivalued metric space (spatial metric space) is introduced and this space reduces to  metric space in the real domain. 
\end{Remark}


\begin{Theorem}\label{3.4}
	The spatial distance measure $(D_1,D_2,D_3)$ is normable if $D_1,D_2$, and $D_3$ satisfies translation invariant and homogeneity conditions.
\end{Theorem}
\begin{Proof}
	Let us show that $D_1$ is translation invariant.
	
	Let $\tilde{{\underline{y}}}=(\mu_{j},\nu_{j})$, $\tilde{{\underline{x}}}=(\mu_{j}^{1},\nu_{j}^{1})$, and $\tilde{{\underline{z}}}=(\mu_{j}^{2},\nu_{j}^{2})$
	\begin{itemize}
		\item[1.] To show, $D_1(\tilde{{\underline{y}}}+\tilde{{\underline{x}}}, 	\tilde{{\underline{x}}}+\tilde{{\underline{z}}})= 	D_1(\tilde{{\underline{y}}},\tilde{{\underline{z}}})$
		
		Now,~ $\tilde{{\underline{y}}}+\tilde{{\underline{x}}}=(\mu_{j}+\mu_{j}^{1},\nu_{j}+\nu_{j}^{1})$, ~and~~$\tilde{{\underline{x}}}+\tilde{{\underline{z}}}=(\mu_{j}^{1}+\mu_{j}^{2},\nu_{j}^{1}+\nu_{j}^{2})$\\ 
		$\mbox{So},~~D_1(\tilde{{\underline{y}}}+\tilde{{\underline{x}}}, 	\tilde{{\underline{x}}}+\tilde{{\underline{z}}})$
		\begin{eqnarray}
			&=&\frac{1}{4k}\sum\limits_{j=1}^{k}\bigg(\big|\Delta_{11}^{1}\big((\mu_j+\mu_j^{1},0),(\mu_{j-1}+\mu_{j-1}^{1},0)\big)
			-\Delta_{11}^{1}\big((\mu_j^{1}+\mu_{j}^{2},0),(\mu_{j-1}^{1}+\mu_{j-1}^{2},0)\big)\big| \nonumber \\
			&&\hspace{0.0cm}+\big|\Delta_{11}^{1}\big((\mu_j+\mu_{j}^{1},\nu_{j}+\nu_{j}^{1}),(\mu_{j-1}+\mu_{j-1}^{1},\nu_{j-1}+\nu_{j-1}^{1})\big)
			-\Delta_{11}^{1}\big((\mu_j^{1}+\mu_{j}^{2},\nu_j^{1}+\nu_{j}^{2}),(\mu_{j-1}^{1}+\mu_{j-1}^{2},\nu_{j-1}^{1}+\nu_{j-1}^{2})\big)\big|\bigg) \nonumber\\
	&=&\frac{1}{4k}\sum\limits_{j=1}^{k}\bigg(\big|\big(\mu_j+\mu_j^{1}-\mu_{j-1}-\mu_{j-1}^{1}+0-0\big)
			-\big(\mu_j^{1}+\mu_{j}^{2}-\mu_{j-1}^{1}-\mu_{j-1}^{2}-0+0\big)\big| \nonumber \\
			&&+\big|\big(\mu_j+\mu_{j}^{1}-\mu_{j-1}-\mu_{j-1}^{1}+\nu_{j-1}+\nu_{j-1}^{1}-\nu_{j}-\nu_{j}^{1}\big) -\big(\mu_j^{1}+\mu_{j}^{2}-\mu_{j-1}^{1}-\mu_{j-1}^{2}+\nu_{j-1}^{1}+\nu_{j-1}^{2}-\nu_j^{1}-\nu_{j}^{2}\big)\big|\bigg) \nonumber\\
			&=&\frac{1}{4k}\sum\limits_{j=1}^{k}\bigg(\big|\big(\mu_j-\mu_{j-1}\big)
			-\big(\mu_{j}^{2}-\mu_{j-1}^{2}\big)\big|+\big|\big(\mu_j-\mu_{j-1}+\nu_{j-1}-\nu_{j}\big)
			-\big(\mu_{j}^{2}-\mu_{j-1}^{2}+\nu_{j-1}^{2}-\nu_{j}^{2}\big)\big|\bigg) \nonumber\\
			&=&\frac{1}{4k}\sum\limits_{j=1}^{k}\bigg(\big|\Delta_{11}^{1}\big((\mu_j,0),(\mu_{j-1},0)\big)
			-\Delta_{11}^{1}\big((\mu_j^{1},0),(\mu_{j-1}^{1},0)\big)\big| \nonumber \\
			&&\hspace{2.5cm}+\big|\Delta_{11}^{1}\big((\mu_j,\nu_{j}),(\mu_{j-1},\nu_{j-1})\big)
			-\Delta_{11}^{1}\big((\mu_j^{2},\nu_j^{2}),(\mu_{j-1}^{2},\nu_{j-1}^{2})\big)\big|\bigg) \nonumber\\
			&=&D_1(\tilde{{\underline{y}}},\tilde{{\underline{z}}}) \nonumber
		\end{eqnarray}
		
		\item[2.] To show, $D_1(\eta\tilde{{\underline{y}}},\eta\tilde{{\underline{x}}})=|\eta| D_1(\tilde{{\underline{y}}},\tilde{{\underline{x}}})$
		\begin{eqnarray}
			D_1(\eta\tilde{{\underline{y}}},\eta\tilde{{\underline{x}}})&=&\frac{1}{4k}\sum\limits_{j=1}^{k}\bigg(\big|\Delta_{11}^{1}\big(
			(\eta \mu_j,0),(\eta \mu_{j-1},0)\big)
			-\Delta_{11}^{1}\big((\eta \mu_j^{1},0),(\eta \mu_{j-1}^{1},0)\big)\big| \nonumber \\
			&&\hspace*{2.5cm}+\big|\Delta_{11}^{1}\big((\eta \mu_j,\eta \nu_{j}),(\eta \mu_{j-1},\eta \nu_{j-1})\big)-\Delta_{11}^{1}\big((\eta \mu_j^{1},\eta \nu_j^{1}),(\eta \mu_{j-1}^{1},\eta \nu_{j-1}^{1})\big)\big|\bigg) \nonumber\\
			&=&\frac{1}{4k}\sum\limits_{j=1}^{k}\bigg(\big|\big(\eta \mu_j-\eta \mu_{j-1}+0-0\big)
			-\big(\eta \mu_j^{1}-\eta \mu_{j-1}^{1}-0+0\big)\big| \nonumber \\
			&&\hspace{2.5cm}+\big|\big(\eta \mu_j-\eta \mu_{j-1}+\eta \nu_{j-1}-\eta \nu_{j}\big)
			-\big(\eta \mu_j^{1}-\eta \mu_{j-1}^{1}+\eta \nu_{j-1}^{1}-\eta \nu_j^{1}\big)\big|\bigg) \nonumber\\
			&=&\frac{1}{4k}\sum\limits_{j=1}^{k}\bigg(\big|\eta\big(\mu_j-\mu_{j-1}+0-0\big)
			-\eta\big(\mu_j^{1}-\mu_{j-1}^{1}-0+0\big)\big| \nonumber \\
			&&\hspace{2.5cm}+\big|\eta\big(\mu_j-\mu_{j-1}+\nu_{j-1}-\nu_{j}\big)
			-\eta\big(\mu_j^{1}-\mu_{j-1}^{1}+\nu_{j-1}^{1}-\nu_j^{1}\big)\big|\bigg) \nonumber\\
			&=&|\eta|\bigg(\frac{1}{4k}\sum\limits_{j=1}^{k}\bigg(\big|\big(\mu_j-\mu_{j-1}+0-0\big)
			-\big(\mu_j^{1}-\mu_{j-1}^{1}-0+0\big)\big| \nonumber \\
			&&\hspace{2.5cm}+\big|\big(\mu_j-\mu_{j-1}+\nu_{j-1}-\nu_{j}\big)
			-\big(\mu_j^{1}-\mu_{j-1}^{1}+\nu_{j-1}^{1}-\nu_j^{1}\big)\big|\bigg)\bigg) \nonumber\\
			&=&|\eta|D_1(\tilde{{\underline{y}}},\tilde{{\underline{x}}}) \nonumber
		\end{eqnarray}
	\end{itemize}
	Hence, $D_1$ is normable. On the similar lines $D_2$ and $D_3$ are also normable. Thus, the spatial distance measure $(D_1,D_2,D_3)$ is normable. $~~~~~~~~~~~~~\square$
\end{Proof}

\begin{Remark}
	Theorem-\ref{3.4} gives an existence to following three norms: 
	\begin{eqnarray}
		\|\tilde{{\underline{y}}}\|_{1}&=&	 \frac{1}{4k}{\sum\limits_{j=1}^{k}\bigg(|\Delta_{11}^{1}\big((\mu_j,0),(\mu_{j-1},0)\big)|}+\big|\Delta_{11}^{1}\big((\mu_j,\nu_j),(\mu_{j-1},\nu_{j-1})\big)\big|\bigg)\\
		\|\tilde{{\underline{y}}}\|_{2}&=&	 \frac{1}{4k}{\sum\limits_{j=1}^{k}\bigg(|\Delta_{11}^{1}\big((0,\nu_j),(0,\nu_{j-1})\big)\big|}+\big|\Delta_{11}^{1}\big((\mu_j,\nu_j),(\mu_{j-1},\nu_{j-1})\big)\big|\bigg)\\
		\|\tilde{{\underline{y}}}\|_{3}&=&	 \frac{1}{4k}{\sum\limits_{j=1}^{k}\bigg(|\Delta_{11}^{1}\big((\mu_j,0),(\mu_{j-1},0)\big)|+|\Delta_{11}^{1}\big((0,\nu_j),(0,\nu_{j-1})\big)|}\bigg)
	\end{eqnarray}
	Here, $\mu_0=\nu_0=0$.
\end{Remark}

\begin{Theorem} \label{3.5}The spatial normed spaces $\bigg( \underline{y}^k,(\|\tilde{{\underline{y}}}\|_{1},\|\tilde{{\underline{y}}}\|_{2},\|\tilde{{\underline{y}}}\|_{3})\bigg)$ are complete under the membership dominant norm $\|\tilde{{\underline{y}}}\|_{1}$, non-membership dominant norm $\|\tilde{{\underline{y}}}\|_{2}$, and equidominant norm $\|\tilde{{\underline{y}}}\|_{3}$. 
\end{Theorem}

\begin{Proof} The spatial normed space consists of three normed spaces
	$\bigg( \underline{y}^k,\|\tilde{{\underline{y}}}\|_{1}\bigg)$, $\bigg( \underline{y}^k,\|\tilde{{\underline{y}}}\|_{2}\bigg)$, and $\bigg( \underline{y}^k,\|\tilde{{\underline{y}}}\|_{3}\bigg)$. To show that spatial normed space is complete, we must prove all three normed spaces are complete. Let us prove the completeness of $\bigg( \underline{y}^k,\|\tilde{{\underline{y}}}\|_{1}\bigg)$.\\
	Let $\big((\mu_{j}^{e},\nu_{j}^{e})_{j=1}^{\infty}:\mu_{j}^{e},\nu_{j}^{e}\geq 0~\mbox{for}~ 1\leq j\leq k,~\mbox{and}~\mu_{j}^{e},\nu_{j}^{e}=0~\mbox{for}~j>k)\big)_{e=1}^{\infty}$  
	be a Cauchy sequence in $\underline{y}^k$. 
		  Now, for a given $\varepsilon>0$, there exists a positive integer $k_0(\varepsilon)>0$ such that 
		\begin{eqnarray}
			\|(\mu_{j}^{e},\nu_{j}^{e})_{j=1}^{\infty}-(\mu_{j}^{f},\nu_{j}^{f})_{j=1}^{\infty}\|_{1}&=&\frac{1}{4k}\sum\limits_{j=1}^{k}\bigg(\big|\Delta_{11}^{1}\big(\big((\mu^{e}_j,0),(\mu^{e}_{j-1},0)\big)-\big((\mu^{f}_j,0),(\mu^{f}_{j-1},0)\big)\big)\big| \nonumber\\
			&&\hspace*{0.5cm}+
			\big|\Delta_{11}^{1}\big(\big((\mu^{e}_j,\nu^{e}_j),(\mu^{e}_{j-1},\nu^{e}_{j-1})\big)-\big((\mu^{f}_j,\nu^{f}_j),(\mu^{f}_{j-1},\nu^{f}_{j-1})\big)\big)\big|\bigg)<\varepsilon \nonumber
		\end{eqnarray}
		for all $e,f>k_{0}(\varepsilon)$. Therefore
		\begin{eqnarray}
			&&\frac{1}{4k}\bigg(\sum\limits_{j=1}^{k}\bigg|\big|\Delta_{11}^{1}\big((\mu^{e}_j,0),(\mu^{e}_{j-1},0)\big)\big|-\big|\Delta_{11}^{1}\big((\mu^{f}_j,0),(\mu^{f}_{j-1},0)\big)\big|\bigg|+\sum\limits_{j=1}^{k}\bigg|
			\big|\Delta_{11}^{1}\big((\mu^{e}_j,\nu^{e}_j),(\mu^{e}_{j-1},\nu^{e}_{j-1})\big)\big| \nonumber\\
			&&\hspace*{2.5cm}-\big|\Delta_{11}^{1}\big((\mu^{f}_j,\nu^{f}_j),(\mu^{f}_{j-1},\nu^{f}_{j-1})\big)\big|\bigg|\bigg)<\varepsilon ~~ \mbox{for all e,~f }> k_0(\varepsilon). \nonumber
		\end{eqnarray}
		This means $\bigg(\sum\limits_{j=1}^{k}|\Delta_{11}^{1}\big((\mu^{e}_j,0),(\mu^{e}_{j-1},0)\big)|\bigg)_{e=1}^{\infty}$ and $\bigg(\sum\limits_{j=1}^{k}|\Delta_{11}^{1}\big((\mu^{e}_j,\nu^{e}_j),(\mu^{e}_{j-1},\nu^{e}_{j-1})\big)|_{e=1}^{\infty}\bigg)_{e=1}^{\infty}$
		are cauchy sequences in $\mathbb R$. Since $\mathbb R$ is complete, the sequences converge, i.e. 
		$$\sum\limits_{j=1}^{k}\big|\Delta_{11}^{1}\big((\mu^{e}_j,0),(\mu^{e}_{j-1},0)\big)|\rightarrow \sum\limits_{j=1}^{k}\big|\Delta_{11}^{1}\big((\mu_{j},0),(\mu_{j-1},0)\big)|~~~~\mbox{as}~~~~e\rightarrow\infty$$
		and
		$$\sum\limits_{j=1}^{k}|\Delta_{11}^{1}\big((\mu^{e}_j,\nu^{e}_j),(\mu^{e}_{j-1},\nu^{e}_{j-1})\big)|\rightarrow \sum\limits_{j=1}^{k}\big|\Delta_{11}^{1}\big((\mu_{j},\nu_{j}),(\mu_{j-1},\nu_{j-1})\big)|~~~~\mbox{as}~~~~e\rightarrow\infty$$
		As absolute convergence implies convergence in $\mathbb R$, so we have
		$$\sum\limits_{j=1}^{k}\Delta_{11}^{1}\big((\mu^{e}_j,0),(\mu^{e}_{j-1},0)\big)\rightarrow \sum\limits_{j=1}^{k}\Delta_{11}^{1}\big((\mu_{j},0),(\mu_{j-1},0)\big)~~~~\mbox{as}~~~~e\rightarrow\infty$$
		and
		$$\sum\limits_{j=1}^{k}\Delta_{11}^{1}\big((\mu^{e}_j,\nu^{e}_j),(\mu^{e}_{j-1},\nu^{e}_{j-1})\big)\rightarrow \sum\limits_{j=1}^{k}\Delta_{11}^{1}\big((\mu_{j},\nu_{j}),(\mu_{j-1},\nu_{j-1})\big)~~~~\mbox{as}~~~~e\rightarrow\infty$$
		Hence we have,
		$$ \lim\limits_{e\rightarrow\infty}\|(\mu_{j}^{e},\nu_{j}^{e})_{j=1}^{\infty}-(\mu_{j},\nu_{j})_{j=1}^{\infty}\|_{1}=0$$
		
		Let $x^{e} =\sum\limits_{j=1}^{k}|\Delta_{11}^{1}\big((\mu^{e}_j,0),(\mu^{e}_{j-1},0)\big)|$, $z^{e} =\sum\limits_{j=1}^{k}|\Delta_{11}^{1}\big((\mu^{e}_j,\nu^{e}_j),(\mu^{e}_{j-1},\nu^{e}_{j-1})\big)|$ and so $x^{e}+z^{e}$ is an element of ${e}_{\infty}$. Therefore
		\begin{eqnarray}
			sup_{e \rightarrow \infty}~~
			\frac{1}{4k}{\sum\limits_{j=1}^{k}\bigg(|\Delta_{11}^{1}\big((\mu_j^e,0),(\mu_{j-1}^e,0)\big)|}+\big|\Delta_{11}^{1}\big((\mu_j^e,\nu_j^e),(\mu_{j-1}^e,\nu_{j-1}^e)\big)\big|\bigg)	
			 ~\leq k_1. \nonumber
		\end{eqnarray}
	\begin{eqnarray}
\mbox{Since,~~~}\|\big(\mu_{j}^{},\nu_{j}^{}\big)_{j=1}^{\infty}\|_{1}	\leq\|\big(\mu_{j}^{},\nu_{j}^{}\big)_{j=1}^{\infty}-\big(\mu_{j}^{e},\nu_{j}^{e}\big)_{j=1}^{\infty}\|_{1}+\|\big(\mu_{j}^{e},\nu_{j}^{e}\big)_{j=1}^{\infty}\|_{1}\nonumber
		~~\leq \varepsilon  +k_1.
	\end{eqnarray}
		It follows that $\bigg((\mu_{j}^{e},\nu_{j}^{e})_{j=1}^{\infty}\bigg)_{e=1}^{\infty}\in{\underline{y}^k}$. Since $\bigg((\mu_{j}^{e},\nu_{j}^{e})_{j=1}^{\infty}\bigg)_{e=1}^{\infty}$ was an arbitrary Cauchy sequence, the space $\bigg( \underline{y}^k,\|\tilde{{\underline{y}}}\|_{1}\bigg)$ is complete. Similarly we can establish the completeness of $\bigg( \underline{y}^k,\|\tilde{{\underline{y}}}\|_{2}\bigg)$ and $\bigg( \underline{y}^k,\|\tilde{{\underline{y}}}\|_{3}\bigg)$. The completeness of spatial normed space  $\bigg( \underline{y}^k,(\|\tilde{{\underline{y}}}\|_{1},\|\tilde{{\underline{y}}}\|_{2},\|\tilde{{\underline{y}}}\|_{3})\bigg)$ is complete. $~~~~~~~~~~~~~\square$
	\end{Proof}

\begin{definition}
In AIFS, distance measure induces a similarity measure due to the existence of 
complementary relationship. The three distance measures of spatial metric space are used to propose a three branched spatial similarity measure (SSM). The three branches are (1) membership dominant spatial similarity measure $(\mbox{MD-}S_{\rm bv_{}})$, (2) non-membership dominant spatial similarity measure $(\mbox{NMD-}S_{\rm bv_{}})$, (3) equidominant spatial similarity measure $(\mbox{ED-}S_{\rm bv_{}})$.
The proposed SSM has following three branches:
\begin{eqnarray}
	&\mbox{SSM}=& \left\{
	\begin{array}{lr}
		\mbox{MD-}S_{\rm bv_{}}, ~ &  \text{membership dominant  } \\
		\mbox{NMD-}S_{\rm bv_{}},~ & \text{non-membership dominant  }  \\
		\mbox{ED-}S_{\rm bv_{}},~ & \text{equidominant } 
	\end{array}
	\right.
\end{eqnarray}
where,
		\begin{eqnarray}
     \mbox{MD-}S_{\rm bv_{}}&=& 1-	 \frac{1}{4k}\sum\limits_{j=1}^{k}\bigg(\big|\Delta_{11}^{1}\big((\mu_j,0),(\mu_{j-1},0)\big)
	-\Delta_{11}^{1}\big((\mu_j^{1},0),(\mu_{j-1}^{1},0)\big)\big| \nonumber \\
	&&\hspace*{2.5cm}+\big|\Delta_{11}^{1}\big((\mu_j,\nu_j),(\mu_{j-1},\nu_{j-1})\big)
	-\Delta_{11}^{1}\big((\mu_j^{1},\nu_j^{1}),(\mu_{j-1}^{1},\nu_{j-1}^{1})\big)\big|\bigg)\hspace*{0.6cm} \nonumber \\
	\mbox{NMD-}S_{\rm bv_{}}&=& 1-	 \frac{1}{4k}\sum\limits_{j=1}^{k}\bigg(\big|\Delta_{11}^{1}\big((0,\nu_j),(0,\nu_{j-1}
	)\big)
	-\Delta_{11}^{1}\big((0,\nu_j^{1}),(0,\nu_{j-1}^{1})\big)\big| \nonumber\\
	&&\hspace*{2.5cm}+\big|\Delta_{11}^{1}\big((\mu_j,\nu_j),(\mu_{j-1},\nu_{j-1})\big)
	-\Delta_{11}^{1}\big((\mu_j^{1},\nu_j^{1}),(\mu_{j-1}^{1},\nu_{j-1}^{1})\big)\big|\bigg)	\hspace*{0.6cm} \nonumber \\
	\mbox{ED-}S_{\rm bv_{}}&=& 1-\frac{1}{4k}\sum\limits_{j=1}^{k}\bigg(\big|\Delta_{11}^{1}\big((\mu_j,0),(\mu_{j-1},0)\big)
	-\Delta_{11}^{1}\big((\mu_j^{1},0),(\mu_{j-1}^{1},0)\big) \big|  \nonumber\\
	&&\hspace{2.5cm}+\big|\Delta_{11}^{1}\big((0,\nu_j),(0,\nu_{j-1})
	\big)-\Delta_{11}^{1}\big((0,\nu_j^{1}),(0,\nu_{j-1}^{1})\big)\big| \bigg)\nonumber
\end{eqnarray}
\end{definition}
\begin{Remark}
	The proposal of SSM does not involve weights as it is induced from spatial metric space.\\ In $\bigg(\underline{y}^{k},(D_1,D_2,D_3)\bigg)$, the weights are intrinsically added. Moreover, $\bigg(\underline{y}^{k},\rm SSM\bigg)$ is also complete.
\end{Remark}


\section{Experimental simulations and comparisons}\label{sec4}
Generally, decision-making, pattern recognition, clustering, and classification are the domains that are used to check the proficiencies of the similarity/distance measures (see: \cite{57},\cite{58},\cite{59},\cite{60},\cite{61}).
Experimental validation of the proposed spatial similarity measure is conducted over some problems of pattern recognition for its comprehensive analysis. Here, we have coined a term called strong classification to clarify the findings of SSM.
	\begin{definition} Strong Classification \\
		If the classification of the pattern obtained using $\mbox{ED-}S_{\rm bv_{}}$ branch is verified by $\mbox{MD-}S_{\rm bv_{}}$ and $\mbox{NMD-}S_{\rm bv_{}}$ branches, then such a classification is referred as strong classification.   
\end{definition}
 

\subsection{Pattern recognition problems}

The pattern recognition problem identifies an unknown AIFS pattern while calculating its similarity values with the known AIFSs. The highest similarity value is used for the characterization of unknown AIFS. Here, four benchmark numeric data examples and two medical diagnosis problems of pattern recognition are dealt with. We have used the proposed SSM to categorize the unknown AIFS in all four benchmark numeric data examples (\cite{10}, \cite{24}, \cite{52}). Further, we have shown that SSM properly works over real-life medical problems, namely, medical diagnosis and cancer diagnosis (\cite{2}, \cite{23},\cite{24},\cite{55}). Here, a thorough experimental validation of SSM is also carried out. We have considered the highest similarity value to recognize the unknown pattern. \\  
~\\
{\bf{Example 1.}} (\cite{2}, \cite{8}, \cite{10}, \cite{12}, \cite{24}): In a universe of discourse $U = \{y_1, y_2, y_3,y_4\}$ the three patterns are in the form of AIFS $B_1$, $B_2$, $B_3$ as follows:   
\begin{center}
	\scriptsize
	\centering
	$B_1=\{ (y_1,0.5,0.3), (y_2,0.7,0.0), (y_3,0.4,0.5), (y_4,0.7,0.3)\}$\\
	$B_2=\{(y_1,0.5,0.2), (y_2,0.6,0.1), (y_3,0.2,0.7), (y_4,0.7,0.3)       \}$\\
	$B_3=\{(y_1,0.5,0.4), (y_2,0.7,0.1), (y_3,0.4,0.6), (y_4,0.7,0.2)\}$
\end{center}
We have to adjudge the classification of the unknown AIFS 
\begin{table}[H]
	\scriptsize
	\centering
	$A=\{ (y_1,0.4,0.3), (y_2,0.7,0.1), (y_3,0.3,0.6), (y_4,0.7,0.3) \}.$
\end{table}
\hspace*{-0.5cm}\textbf{Result Discussion:} Let us calculate the similarity between $B_1$ and $A$, $B_2$ and $A$, $B_3$ and $A$, using SSM. The proposed SSM consists of three branches: membership dominant, non-membership dominant, and equidominant similarity measure. These branches result in three similarity measures. 
The three similarity measures yield three values while calculating the similarity between $B_i$ $(1\leq i \leq 3 )$ and $A$.\\
Example-1 is already addressed in several papers, and the classification results corresponding to some well-known similarity measures are known (see: \cite{12}, \cite{24},\cite{2}). These classification results are comapred with those obtained using $\mbox{MD-}S_{\rm bv_{}}$, $\mbox{NMD-}S_{ \rm bv_{}}$, and  $\mbox{ED-}S_{\rm bv_{}}$ in Table-\ref{Example1}. The $\mbox{ED-}S_{\rm bv_{}}$ is same as $S_{\rm Az_{1}}$ so both are classifying the unknown pattern $A$ into $B_1$, and $S_{\rm Az_{2}}$ is also classifying the same. The two branches $\mbox{MD-}S_{\rm bv_{}}$ and $\mbox{NMD-}S_{\rm bv_{}}$ of proposed spatial similarity measure are classifying the unknown pattern $A$ into $B_3$. Most of the classifying similarity measures such as $S_{\rm Bd}, S_{\rm C}, S_{\rm Fz}, S_{\rm Dc},$ $S_{\rm Hm}$, $S_{\rm Az_{p}^h}$ classify the unknown pattern $A$ into $B_3$.
The similarity does not identify the unknown pattern $A$ measures $S_{\rm L_1}, S_{\rm Hk}, S_{\rm L_2}, S_{\rm Lzd}, S_{\rm M}, S^{1}_{\rm Hy_1}, S^{1}_{\rm Hy_2}, S^{1}_{\rm Hy_3}, S^{2}_{\rm Hy_1}, S^{2}_{\rm Hy_2}, S^{2}_{\rm Hy_3}, S^{3}_{\rm Hy},$ $S_{\rm Ls}$ as their maximal similarity value is equal for two or three given patterns (see: Table-\ref{Example1}). Our SSM is also efficient as it accurately characterizes the unknown pattern by exploiting the dominating behavior of either membership or non-membership.\\ 

\begin{table}[htp]
	\centering
	\fontsize{7}{7}\selectfont
		\caption{Example 1: The classification obtained by the proposed SSM and other well-known measures}
		\renewcommand{\arraystretch}{1.5}	\label{Example1}
	\begin{tabular}{|p{2.5cm}|p{1.5cm}|p{1.5cm}|p{1.5cm}|p{2.0cm}|}
		\hline 
		Similarity Measures &  $S(B_1,A)$ & $S(B_2,A)$ & $S(B_3,A)$& Result\\
		\hline
	$ S_{\rm L_{1}}$ & 0.9500 & 0.9375 & 0.9500 & Unclassified\\
		\hline
	$ S_{\rm Hk}$ & 0.9500 & 0.9375 & 0.9500 & Unclassified\\
		\hline
	$ S_{\rm L_{2}}$ & 0.9293 & 0.9209 & 0.9293 &Unclassified\\
		\hline
	$ S_{\rm Lzd}$ &  0.9293 & 0.9209 & 0.9293 &Unclassified\\
		\hline
	$ S_{\rm M}$~~($p=1$) & 0.9500 & 0.9375 & 0.9500 &Unclassified\\
		 \hline
    $ S_{\rm Fz}$ &  0.9500 & 0.9375 & 0.9625 &$~~~~B_3$\\
		 \hline
	 $ S_{\rm C}$  & 0.9500 & 0.9375 & 0.9750 &$~~~~B_3$\\
	     \hline
	    $ S^{1}_{\rm Hy_1}$ & 0.9250 & 0.9250 & 0.9250 &Unclassified\\
	      \hline
		$ S^{1}_{\rm Hy_2}$  &0.8857  &0.8857  &0.8857  &Unclassified\\
	\hline
	$ S^{1}_{\rm Hy_3}$ &0.8605  & 0.8605 & 0.8605 &Unclassified\\
	\hline
    $ S^{2}_{\rm Hy_1}$~~($p=1$) & 0.9500 & 0.9375 & 0.9500 &Unclassified\\
     \hline
    $ S^{2}_{\rm Hy_2}$~~($p=1$) &0.8899  &0.8641  & 0.8899 &Unclassified\\
    \hline
    $ S^{2}_{\rm Hy_3}$~~($p=1$)  &0.8636  & 0.8333 &0.8636  &Unclassified\\
    \hline
    $ S^{3}_{\rm Hy}$ &0.9500  &0.9375  &0.9500  &Unclassified\\
    \hline
	$ S_{\rm Bd}$~~($p=1,z=2$)  &0.9500  & 0.9375 &0.9667  &$~~~~B_3$\\
      \hline
     $ S_{\rm Dc}$~~($p=1$) &0.9500  &0.9375  &0.9750  &$~~~~B_3$\\
     \hline
     $ S_{\rm Ls}$~~($p=1$) &0.9500  &0.9375  &0.9500  &Unclassified\\
     \hline
     $ S_{\rm Hm}$ &0.9500  &0.9333  &0.9583  &$~~~~B_3$\\
     \hline
      $S_{\rm Az_{1}}$ &0.9250 &0.9000 &0.9125 &$~~~~B_1$ \\
     \hline
     $S_{\rm Az_{2}}$ &0.9134 &0.8882 &0.9065 & $~~~~B_1$ \\
     \hline
     $S_{\rm Az_{1}}^{h}$ &0.9250 &0.9125 &0.9375 & $~~~~B_3$ \\
     \hline
     $S_{\rm Az_{2}}^{h}$ &0.9134 &0.9065 &0.9209 &$~~~~B_3$  \\
     \hline
     $\mbox{MD-}S_{\rm bv_{}}$ &0.9500 &0.9250 & 0.9688&$~~~~B_3$ \\
     \hline
      $\mbox{NMD-}S_{\rm bv_{}}$ &0.9625 &0.9250 &0.9750 & $~~~~B_3$\\
     \hline
      $\mbox{ED-}S_{\rm bv_{}}$ &0.9625 &0.9500 &0.9563 &$~~~~B_1$ \\
     \hline
 \end{tabular}
\end{table}
\hspace*{-0.5cm}{\bf{Example 2.}} (\cite{10},\cite{24},\cite{36},\cite{35}): The three AIFS patterns namely $B_1$, $B_2$ and $B_3$ are given in the universe of discourse $U=\{y_1,y_2,y_3\}$ 
\begin{table}[H]
	\scriptsize
	\centering
	$B_1=\{ (y_1,0.34,0.34), (y_2,0.19,0.48), (y_3,0.02,0.12)\}$
	~\\
	$B_2=\{(y_1,0.35,0.33), (y_2,0.20,0.47), (y_3,0.00,0.14)\}$
	~\\
	$B_3=\{(y_1,0.33,0.35), (y_2,0.21,0.46), (y_3,0.01,0.13)\}$
\end{table}

We have to classify the unknown pattern $A$, where
\begin{table}[H]
	\scriptsize
	\centering
	$A=\{ (y_1,0.37,0.31), (y_2,0.23,0.44), (y_3,0.04,0.10) \}.$
\end{table}
\hspace*{-0.5cm}\textbf{Result Discussion:} The aim is to classify the unknown pattern $A$ with the help of known patterns $B_1,~ B_2$ and $ B_3$. Let us calculate the similarities between known and unknown patterns using the proposed SSM. The three branches $\mbox{MD-}S_{\rm bv_{}}$, $\mbox{NMD-}S_{\rm bv_{}}$, and  $\mbox{ED-}S_{\rm bv_{}}$ of the SSM are classifying the the unknown pattern $A$ into $B_2$ (see: Table-\ref{Example2}).
It is observed (see: \cite{10},\cite{24},\cite{36},\cite{2}) that most of the similarity measures fail to identify the pattern $A$. Here, one of the similarity measure $S_{\rm Az_{p}^h}$ \cite{2} is classifying the unknown pattern $A$ into $B_1$. The SSM is strongly classifying $A$ into $B_2$ due to the similar classification from all the branches (see: Table-\ref{Example2}). The strong classification means the obtained result is highly authentic. 
\begin{table}[htp]
	\centering
	\fontsize{7}{7}\selectfont
	\renewcommand{\arraystretch}{1.5}
	\caption{Example 2: The classification obtained by the proposed SSM and other well-known measures}
		\renewcommand{\arraystretch}{1.5}\label{Example2}
	\begin{tabular}{|p{2.5cm}|p{1.5cm}|p{1.5cm}|p{1.5cm}|p{2.0cm}|}
		\hline 
		Similarity Measures &  $S(B_1,A)$ & $S(B_2,A)$ & $S(B_3,A)$&~~ Result\\
		\hline
		$ S_{\rm L_{1}}$ & 0.9700 & 0.9700 & 0.9700 & Unclassified \\
		\hline
	    $ S_{\rm Hk}$ & 0.9700 & 0.9700 & 0.9700 & Unclassified\\
	\hline
		$ S_{\rm L_{2}}$ & 0.9689 & 0.9689 & 0.9689 & Unclassified\\
		\hline
		$ S_{\rm Lzd}$& 0.9689 & 0.9689 & 0.9689 & Unclassified\\
		\hline
		$ S_{\rm M}$~~($p=1$)  & 0.9700 & 0.9700 & 0.9700 & Unclassified\\
		\hline
		$ S_{\rm Fz}$ & 0.9700 & 0.9700 & 0.9700 & Unclassified\\
		\hline
		$ S_{\rm C}$  & 0.9700 & 0.9700 & 0.9700 &Unclassified\\
		\hline
		$ S^{1}_{\rm Hy_1}$&0.9700  &0.9700  &0.9700  &Unclassified\\
		\hline
		$ S^{1}_{\rm Hy_2}$ &0.9532  &0.9532  &0.9532  &Unclassified\\
		\hline
		$ S^{1}_{\rm Hy_3}$&0.9417  &0.9417  &0.9417  &Unclassified\\
		\hline
		$ S^{2}_{\rm Hy_1}$~~($p=1$) & 0.9700 & 0.9700 & 0.9700 &Unclassified\\
		\hline
		$ S^{2}_{\rm Hy_2}$~~($p=1$)  &0.9326  &0.9326 & 0.9326 &Unclassified\\
		\hline
		$ S^{2}_{\rm Hy_3}$~~($p=1$) &0.9151  & 0.9151 &0.9151  &Unclassified\\
		\hline
		$ S^{3}_{\rm Hy}$  &0.9700  &0.9700  &0.9700  &Unclassified\\
		\hline
		$ S_{\rm Bd}$~~($p=1,z=2$)&0.9700  &0.9700  &0.9700  &Unclassified\\
		\hline
		$ S_{\rm Dc}$~~($p=1$)  &0.9700  &0.9700  & 0.9700 &Unclassified\\
		\hline
		$ S_{\rm Ls}$~~($p=1$) &0.9700  &0.9700  &0.9700  &Unclassified\\
		\hline
		$ S_{\rm Hm}$ &0.9700  &0.9700  &0.9700  &Unclassified\\
		\hline
		 $S_{\rm Az_{1}}$ &0.9800 &0.9867 &0.9767 &$~~~~B_2$ \\
		\hline
		$S_{\rm Az_{2}}$ &0.9784 &0.9859 &0.9742 &$~~~~B_2$  \\
		\hline
		$S_{\rm Az_{1}}^{h}$ &0.9800 &0.9767 &0.9700 & $~~~~B_1$ \\
		\hline
		$S_{\rm Az_{2}}^{h}$ &0.9784 &0.9735 &0.9689 &$~~~~B_1$  \\
		\hline
		$\mbox{MD-}S_{\rm bv_{}}$ &0.9850 &0.9900 &0.9825 &$~~~~B_2$ \\
		\hline
		$\mbox{NMD-}S_{\rm bv_{}}$ &0.9850 & 0.9900&0.9825 & $~~~~B_2$\\
		\hline
		$\mbox{ED-}S_{\rm bv_{}}$ &0.9900 &0.9933 &0.9883 &$~~~~B_2$ \\
		\hline
	\end{tabular}
\end{table}
~\\
{\bf{Example 3.}} (\cite{14},\cite{24},\cite{36},\cite{53}): 
 Three known patterns, $B_1$, $B_2$, and $B_3$, are represented in terms of AIFS. Let $U= \{y_1, y_2, y_3\}$ be the universe of discourse. We have,
 \begin{table}[H]
 	\scriptsize
 	\centering
 	$B_1=\{ (y_1,1.0,0.0), (y_2,0.8,0.0), (y_3,0.7,0.1)\}$
 	~\\
 	$B_2=\{(y_1,0.8,0.1), (y_2,1.0,0.0), (y_3,0.9,0.0)\}$
 	~\\
 	$B_3=\{(y_1,0.6,0.2), (y_2,0.8,0.0), (y_3,1.0,0.0)\}$
 \end{table}
 We need to classify the unknown pattern $A$ based on the patterns $B_1$, $B_2$, and $B_3$. We have,
 \begin{table}[H]
 	\scriptsize
 	\centering
 	$A=\{ (y_1,0.5,0.3), (y_2,0.6,0.2), (y_3,0.8,0.1)\}$
 \end{table}
 \hspace*{-0.5cm}\textbf{Result Discussion:} The objective is to classify the unknown pattern $A$ with the help of known patterns $B_1,~ B_2$ and $ B_3$. Let us calculate similarities between the patterns using the proposed SSM. In Table-\ref{Example3}, the $\mbox{MD-}S_{\rm bv_{}}$, $\mbox{NMD-}S_{\rm bv_{}}$, and $\mbox{ED-}S_{\rm bv_{}}$ are classifying the pattern $A$ into $B_3$. The research work of \cite{2},\cite{24},\cite{53} suggest that almost all of the similarity measures have classified the unknown pattern $A$ into $B_3$. So, the results obtained from the proposed SSM are in accordance with the existing results (see: Table-\ref{Example3}). Thus, we can say that our proposed SSM strongly characterizes the unknown pattern.\\
	\begin{table}[htp]
	\centering
	\fontsize{7}{7}\selectfont
	\caption{Example 3: The classification obtained by the proposed SSM and other well-known measures}
	\renewcommand{\arraystretch}{1.5}\label{Example3}
	\begin{tabular}{|p{2.5cm}|p{1.5cm}|p{1.5cm}|p{1.5cm}|p{2.0cm}|}
		\hline 
		Similarity Measures &  $S(B_1,A)$ & $S(B_2,A)$ & $S(B_3,A)$& Result\\
		\hline
		$ S_{\rm L_{1}}$ & 0.7833 & 0.7833 & 0.8500 &$~~~~B_3$\\
		\hline
		$ S_{\rm Hk}$ & 0.7833 & 0.7833 & 0.8500 &$~~~~B_3$\\
		\hline
		$ S_{\rm L_{2}}$ & 0.7323 & 0.7585 & 0.8419 &$~~~~B_3$\\
		\hline
		$ S_{\rm Lzd}$  & 0.7323 & 0.7585 & 0.8419 &$~~~~B_3$\\
		\hline
		$ S_{\rm M}$~~($p=1$)& 0.7833 & 0.7833 & 0.8500 &$~~~~B_3$\\
		\hline
		$ S_{\rm Fz}$ & 0.7833 & 0.7833 & 0.8500 &$~~~~B_3$\\
		\hline
		$ S_{\rm C}$  & 0.7833 & 0.7833 & 0.8500 &$~~~~B_3$\\
		\hline
			$ S^{1}_{\rm Hy_1}$ &0.7333  &0.7333  &0.8333  &$~~~~B_3$\\
		\hline
		$ S^{1}_{\rm Hy_2}$ &0.6297  &0.6297  &0.7571  &$~~~~B_3$\\
		\hline
		$ S^{1}_{\rm Hy_3}$ &0.5789  &0.5789  & 0.7143 &$~~~~B_3$\\
		\hline
	$ S^{2}_{\rm Hy_1}$~~($p=1$)  & 0.7833 & 0.7833 & 0.8500 &$~~~~B_3$\\
		\hline
		$ S^{2}_{\rm Hy_2}$~~($p=1$) &0.5933  &0.5933  & 0.7003 &$~~~~B_3$\\
		\hline
		$ S^{2}_{\rm Hy_3}$~~($p=1$) &0.5465  & 0.5465 &0.6538  &$~~~~B_3$\\
		\hline
		$ S^{3}_{\rm Hy}$  &0.7833  &0.7833  &0.8500  &$~~~~B_3$\\
		\hline
		$ S_{\rm Bd}$~~($p=1,z=2$)&0.7833  &0.7833  &0.8500  &$~~~~B_3$\\
		\hline
		$ S_{\rm Dc}$~~($p=1$)&0.7833  &0.7833  &0.8500  &$~~~~B_3$\\
		\hline
		$ S_{\rm Ls}$~~($p=1$) & 0.7833 &0.7833  &0.8500  &$~~~~B_3$\\
		\hline
		$ S_{\rm Hm}$ & 0.7667 &   0.7667  &  0.8444 & $~~~~B_3$ \\
		\hline
		 $S_{\rm Az_{1}}$ &0.7167 &0.8333 &0.9167 &$~~~~B_3$ \\
		\hline
		$S_{\rm Az_{2}}$ &0.6918 &0.8000 &0.9087 &$~~~~B_3$  \\
		\hline
		$S_{\rm Az_{1}}^{h}$ &0.8000 &0.8167 &0.9000 &$~~~~B_3$  \\
		\hline
		$S_{\rm Az_{2}}^{h}$ &0.7354 &0.7655 &0.8709 &$~~~~B_3$  \\
		\hline
        $\mbox{MD-}S_{\rm bv_{}}$ &0.7667 &0.8583 &0.9417 &$~~~~B_3$ \\
		\hline
		$\mbox{NMD-}S_{\rm bv_{}}$ &0.8083 &0.8917 &0.9333 &$~~~~B_3$ \\
		\hline
		 $\mbox{ED-}S_{\rm bv_{}}$ &0.8583 &0.9167 &0.9583 &$~~~~B_3$ \\
		\hline
	\end{tabular}
\end{table} 
~\\
{\bf{Example 4.}} (\cite{23},\cite{24}, \cite{2}): Let $B_1$, $B_2$ and $B_3$ be the AIFS based known patterns in the universe of discourse $U = \{y_1, y_2, y_3, y_4, y_5, y_6\}$. We have
\begin{table}[H]
	\scriptsize
	\centering
	$B_1=\{ (y_1,0.94,0.00), (y_2,0.88,0.00), (y_3,0.82,0.00),(y_4,0.78,0.02), (y_5,0.75,0.05), (y_6,0.72,0.08)\}$
	~\\
	$B_2=\{(y_1,0.86,0.07), (y_2,0.92,0.04), (y_3,0.98,0.01), (y_4,0.98,0.00), (y_5,0.95,0.00), (y_6,0.92,0.00)\}$
	~\\
	$B_3=\{(y_1,0.66,0.14), (y_2,0.72,0.08), (y_3,0.78,0.02), (y_4,0.84,0.00), (y_5,0.90,0.00), (y_6,0.96,0.00)\}$
\end{table}
We must classify the unknown pattern $A$ based on the patterns $B_1$, $B_2$ and $B_3$. We have,
\begin{table}[H]
	\scriptsize
	\centering
	$A=\{ (y_1,0.53,0.27), (y_2,0.56,0.24), (y_3,0.59,0.21), (y_4,0.64,0.18), (y_5,0.7,0.15), (y_6,0.76,0.12)\}$
\end{table}
\begin{table}[htp]
	\centering
	\fontsize{7}{7}\selectfont
		\caption{Example 4: The classification obtained by the proposed SSM and other well-known measures}
	\renewcommand{\arraystretch}{1.5}\label{Example4}
	\begin{tabular}{|p{2.5cm}|p{1.5cm}|p{1.5cm}|p{1.5cm}|p{2.0cm}|}
		\hline 
		Similarity Measures & $S(B_1,A)$ & $S(B_2,A)$ & $S(B_3,A)$& Result\\
		\hline
		$ S_{\rm L_{1}}$ & 0.8158 & 0.7600 & 0.8325 &$~~~~B_3$\\
		\hline
		$ S_{\rm Hk}$  & 0.8158 & 0.7600 & 0.8325 &$~~~~B_3$\\
		\hline
		$ S_{\rm L_{2}}$  & 0.7842 & 0.7445 & 0.8301 &$~~~~B_3$\\
		\hline
		$ S_{\rm Lzd}$  & 0.7842 & 0.7445 & 0.8301 &$~~~~B_3$\\
		\hline
		$ S_{\rm M}$~~($p=1$) & 0.8158 & 0.7600 & 0.8325 &$~~~~B_3$\\
		\hline
		$ S_{\rm Fz}$ & 0.8192 & 0.7600 & 0.8325 &$~~~~B_3$\\
		\hline
		$ S_{\rm C}$& 0.8225 & 0.7600 & 0.8325 &$~~~~B_3$\\
		\hline
			$ S^{1}_{\rm Hy_1}$& 0.7900 &0.6950  & 0.8200&$~~~~B_3$\\
		\hline
		$ S^{1}_{\rm Hy_2}$ &0.7003  &0.5841  &0.7394  &$~~~~B_3$\\
		\hline
		$ S^{1}_{\rm Hy_3}$  &0.6529  &0.5326  & 0.6949 &$~~~~B_3$\\
		\hline
		$ S^{2}_{\rm Hy_1}$~~($p=1$)  & 0.8158 & 0.7600 & 0.8325 &$~~~~B_3$\\
		\hline
		$ S^{2}_{\rm Hy_2}$~~($p=1$)  &0.6437  &0.5591  & 0.6708 &$~~~~B_3$\\
		\hline
		$ S^{2}_{\rm Hy_3}$~~($p=1$)  &0.5962  & 0.5135 &0.6236  &$~~~~B_3$\\
		\hline
		$ S^{3}_{\rm Hy}$ &0.8158  &0.7600  &0.8325  &$~~~~B_3$\\
		\hline
		$ S_{\rm Bd}$~~($p=1,z=2$)  &0.8203  &0.7600  &0.8325  &$~~~~B_3$\\
		\hline
		$ S_{\rm Dc}$~~($p=1$)  &0.8225  &0.7600  &0.8325  &$~~~~B_3$\\
		\hline
		$ S_{\rm Ls}$~~($p=1$) &0.8158  &0.7600  &0.8325  &$~~~~B_3$\\
		\hline
		$ S_{\rm Hm}$  &   0.8111  &  0.7383  &  0.8283 &$~~~~B_3$\\
		\hline
		 $S_{\rm Az_{1}}$ &0.8867 &0.9250 &0.9617 &$~~~~B_3$ \\
		\hline
		$S_{\rm Az_{2}}$ &0.8437 &0.8804 &0.9427 &$~~~~B_3$  \\
		\hline
		$S_{\rm Az_{1}}^{h}$ &0.8233 &0.8342 &0.9133 &$~~~~B_3$  \\
		\hline
		$S_{\rm Az_{2}}^{h}$ &0.7894 &0.7895 &0.8848 &$~~~~B_3$  \\
		\hline
		 $\mbox{MD-}S_{\rm bv_{}}$ &0.9075 &0.9367 &0.9733 &$~~~~B_3$ \\
		\hline
		$\mbox{NMD-}S_{\rm bv_{}}$ &0.9225 &0.9508 &0.9708 &$~~~~B_3$ \\
		\hline
		$\mbox{ED-}S_{\rm bv_{}}$ &0.9433 &0.9625 &0.9808 &$~~~~B_3$ \\
		\hline
	\end{tabular}
\end{table} 
\hspace*{-0.5cm}\textbf{Result Discussion:} The objective is to classify the unknown pattern $A$ on the basis of known patterns $B_1,~ B_2$ and $ B_3$. Let us calculate similarities between the patterns using the proposed SSM. In Table-\ref{Example4}, the banches $\mbox{MD-}S_{\rm bv_{}}$, $\mbox{NMD-}S_{\rm bv_{}}$, and  $\mbox{ED-}S_{\rm bv_{}}$ are classifying the pattern $A$ into $B_3$. The research work of \cite{23},\cite{24}, \cite{2} suggest that almost all of the similarity measures have classified the unknown pattern $A$ into $B_3$. So, the results obtained from the SSM are in accordance with the existing results (see: Table-\ref{Example4}). Thus, we can say that our proposed SSM strongly characterizes the unknown pattern.
~\\
{\bf{Example 5.}} (Medical diagnosis: \cite{44},\cite{38},\cite{39},\cite{2}) Medical diagnosis problem is a real-life application of the similarity measure. The application of similarity measures is essential in medical diagnosis, as it involves evaluating the degree of correlation between patients and diseases. This evaluation helps doctors to analyze the patients conditions and enables them to provide appropriate treatments for the respective diseases. 
Let $U$ be the feature space to characterize the nature of the disease. Here,
\begin{center}
	$U=\{y_1=\mbox{Temperature},~y_2=\mbox{Headeche},~y_3=\mbox{Stomach~Pain},~y_4=\mbox{Cough},~y_5=\mbox{Chest~Pain}\}$
\end{center}
Let $B$ be the set of diseases such that:
\begin{center}
	$B= \{B_1=\mbox{Viral~Fever},~B_2=\mbox{Malaria},~B_3=\mbox{Typhoid},~B_4=\mbox{Stomach~Problem},~B_5=\mbox{Heart ~Problem} \}$
\end{center}
which are represented by AIFSs in $U~ (\mbox{universe of discourse}) = \{y_1, y_2, y_3, y_4, y_5\}$ as follows:  
\begin{table}[H]
	\scriptsize
	\centering
	$B_1=\{ (y_1, 0.4,0.0), (y_2, 0.3,0.5), (y_3, 0.1,0.7), (y_4, 0.4,0.3), (y_5, 0.1,0.7)\}$\\
	$B_2=\{(y_1, 0.7,0.0), (y_2, 0.2,0.6), (y_3, 0.0,0.9), (y_4, 0.7,0.0), (y_5, 0.1,0.8)\}$\\
	$B_3=\{(y_1, 0.3,0.3), (y_2, 0.6,0.1), (y_3, 0.2,0.7), (y_4 , 0.2,0.6), (y_5, 0.1,0.9)\}$\\
	$B_4=\{(y_1, 0.1,0.7), (y_2, 0.2,0.4), (y_3, 0.8,0.0), (y_4, 0.2,0.7), (y_5, 0.2,0.7)\}$\\
	$B_5=\{(y_1, 0.1,0.8), (y_2, 0.0,0.8), (y_3, 0.2,0.8), (y_4, 0.2,0.8), (y_5, 0.8,0.1)\}$
\end{table}
A set of patients $A$ is given as follows:
\begin{center}
	$A=\{A_1=\mbox{Al},~A_2=\mbox{Bob},~A_3=\mbox{Joe},~A_4=\mbox{Ted} \}$
\end{center}
which are represented by the AIFSs in $U$, where
\begin{table}[H]
	\scriptsize
	\centering
	$A_1=\{(y_1,0.8,0.1), (y_2,0.6,0.1), (y_3,0.2,0.8), (y_4,0.6,0.1), (y_5,0.1,0.6)    \}$\\
	$A_2=\{(y_1,0.0,0.8), (y_2,0.4,0.4), (y_3,0.6,0.1), (y_4,0.1,0.7), (y_5,0.1,0.8) \}$\\
	$A_3=\{(y_1,0.8,0.1), (y_2,0.8,0.1), (y_3,0.0,0.6), (y_4,0.2,0.7), (y_5,0.0,0.5) \}$\\
	$A_4=\{(y_1,0.6,0.1), (y_2,0.5,0.4), (y_3,0.3,0.4), (y_4,0.7,0.2), (y_5,0.3,0.4) \}$\\
\end{table}
\hspace*{-0.5cm}\textbf{Result Discussion:} The proper diagnosis of the patient $A_i\in A$ $(1\leq i \leq 4)$ depends on the feature space $U$. Let us calculate the similarity degree between the disease $B_j\in B$ $(1\leq j \leq 5)$  and patient $A_i\in A$, $(1\leq i \leq 4)$ using the proposed spatial similarity measure. The branches $\mbox{MD-}S_{\rm bv_{}}$, $\mbox{NMD-}S_{\rm bv_{}}$, and  $\mbox{ED-}S_{\rm bv_{}}$ classify as $A_1 \rightarrow B_1,~ A_2 \rightarrow B_4,~ A_3 \rightarrow B_3,~ A_4 \rightarrow B_1 $ (see: Table-\ref{Example5}). In other words, this classification is strong. Our classification result matches with those obtained by Chen et al. \cite{11}, Ashraf et al. \cite{2}, and Papakostas et al. \cite{39}, and hence the results are justified. A comparison of the obtained result with other similarity measures based on classification is given in Table-\ref{Example5}. 
\begin{table}[htp]
\centering
\fontsize{7}{7}\selectfont
\caption{Example 5: The classification obtained by the proposed SSM and other well-known measures for medical diagnosis problem}
\renewcommand{\arraystretch}{1.5}\label{Example5}
\begin{tabular}{|p{2.5cm}|p{1.5cm}|p{1.5cm}|p{1.5cm}|p{1.5cm}|}
\hline

Similarity Measures  & $A_1$ (Al)& $A_2$ (Bob)& $A_3$ (Joy)& $A_4$ (Ted)\\ \hline

$S_{\rm L_1}$ &$B_2$  &$B_4$   &$B_3$  & $B_1$  \\
 \hline
 
 $S_{\rm Hk}$ &$B_2$  &$B_4$   &$B_3$  & $B_1$  \\
 \hline
 
 $S_{\rm L_2}$ &$B_1$  &$B_4$   &$B_3$  & $B_1$  \\
 \hline
 
  $S_{\rm Lzd}$ &$B_1$  &$B_4$   &$B_3$  & $B_1$  \\
 \hline
 
 $S_{\rm M}$~~($p=1$) &$B_2$  &$B_4$   &$B_3$  & $B_1$  \\
 \hline
 
 $S_{\rm Fz}$   & Unclassified &$B_4$   &$B_3$  & $B_1$  \\
 \hline
 
  $S_{\rm C}$ & $B_1$  &$B_4$   &$B_3$  & $B_1$  \\
 \hline
 
  $S_{\rm Hy_1}^{1}$ & $B_2$  &$B_4$   &$B_3$  & $B_1$  \\
 \hline
 
 $S_{\rm Hy_2}^{1}$ & $B_2$ &$B_4$ &$B_3$  & $B_1$  \\
 \hline
 
 $S_{\rm Hy_3}^{1}$  & $B_2$  &$B_4$   &$B_3$  & $B_1$  \\
 \hline
 
 $S_{\rm Bd}$~~($p=1,z=2$) & $B_1$  &$B_4$  &$B_3$  & $B_1$  \\
 \hline
 
 $S_{\rm Dc}$~~($p=1$)  & $B_1$  &$B_4$ &$B_3$  & $B_1$  \\
 \hline

$S_{\rm Hy_1}^{2}$~~($p=1$) & $B_2$  &$B_4$   &$B_3$  & $B_1$ \\
\hline

$S_{\rm Hy_2}^{2}$~~($p=1$) & $B_2$ &$B_4$ &$B_3$ & $B_1$  \\
\hline

$S_{\rm Hy_3}^{2}$~~($p=1$) & $B_2$ &$B_4$  &$B_3$ & $B_1$  \\
\hline

$S_{\rm Hy}^{3}$ & $B_2$  &$B_4$   &$B_3$  & $B_1$  \\
\hline

$S_{\rm Ls}$~~($p=1$)  &$B_2$  &$B_4$   &$B_3$  & $B_1$  \\
\hline

$S_{\rm Hm}$ & $B_2$  &$B_4$   &$B_3$ & $B_1$  \\
\hline

 $S_{\rm Az_{1}}$ & $B_1$  &$B_4$   &$B_3$  & $B_1$  \\
\hline

$S_{\rm Az_{2}}$ &  $B_1$  &$B_4$   &$B_3$  & $B_1$   \\
\hline

$S_{\rm Az_{1}}^{h}$ &  $B_1$  &$B_4$   &$B_3$  & $B_1$  \\
\hline

$S_{\rm Az_{2}}^{h}$ &  $B_1$  &$B_4$   &$B_3$  & $B_1$   \\
\hline

$\mbox{MD-}S_{\rm bv_{}}$ &$B_1$  &$B_4$  &$B_3$  &$B_1$  \\
\hline
$\mbox{NMD-}S_{\rm bv_{}}$ &$B_1$  &$B_4$  &$B_3$  &$B_1$  \\
\hline
$\mbox{ED-}S_{\rm bv_{}}$ &$B_1$  &$B_4$  &$B_3$  &$B_1$  \\
\hline
	\end{tabular}
\end{table} 

Based on the data presented in Table-\ref{Example5}, it is evident that the proposed SSM along with other similarity measures $S_{\rm L_2}$, $S_{\rm Bd}$, $S_{\rm C}$, $S_{\rm Dc}$, $S_{\rm Lzd}$, $S_{\rm Az_{1}}$, $S_{\rm Az_{2}}$, $S_{\rm Az_{1}}^{h}$, and $S_{\rm Az_{2}}^{h}$ successfully classify the proper diagnosis for the patients as follows: Al $\rightarrow$ Viral Fever, Bob $\rightarrow$ Stomach Problem, Joe $\rightarrow$ Typhoid, and Ted $\rightarrow$ Viral Fever. However, the similarity measures $S_{\rm L_1}$, $S_{\rm Hk}$, $S_{\rm M}$, $S_{\rm Hy_1}^1$, $S_{\rm Hy_2}^1$, $S_{\rm Hy_3}^1$, $S_{\rm Hy_1}^2$, $S_{\rm Hy_2}^2$, $S_{\rm Hy_3}^2$, $S_{\rm Hy}^3$, $S_{\rm Hm}$, and $S_{\rm Ls}$ classify patient $A_1$ for the diagnosis $B_2$, indicating Al $\rightarrow$ Malaria.
Moreover, the similarity measure $S_{\rm Fz}$ encounters issues in predicting the diagnosis for $A_1$ due to division by zero problem.\\

\hspace*{-0.5cm}{\bf{Example 6.}} (Cancer diagnosis: \cite{23},\cite{24},\cite{55},\cite{7}): Let $Y$ be the feature space corresponding to a cancer diagnosis problem, where
\begin{center}
$U=\{y_1=\mbox{Character of stool},~y_2=\mbox{Bellyache},~y_3=\mbox{Ictussileusn},~y_4=\mbox{Chronic sileus},~y_5=\mbox{Anemia}\}$
\end{center}
Let $B=\{B_1=\mbox{Mestastasis},~B_2=\mbox{Recurrence},~B_3=\mbox{Bad},~B_4=\mbox{Well}\}$ be the set of states in which a patient may fall. 
In the $U~\mbox{(universe of discourse)}= \{y_1, y_2, y_3, y_4, y_5\}$, these are given as follows:
\begin{table}[H]
	\scriptsize
	\centering
	$B_1=\{ (y_1, 0.4,0.4), (y_2, 0.3,0.3), (y_3, 0.5,0.1), (y_4, 0.5,0.2), (y_5, 0.6,0.2)\}$\\
	$B_2=\{(y_1, 0.2,0.6), (y_2, 0.3,0.5), (y_3, 0.2,0.3), (y_4, 0.7,0.1), (y_5, 0.8,0.0) \}$\\
	$B_3=\{(y_1, 0.1,0.9), (y_2, 0.0,0.1), (y_3, 0.2,0.7), (y_4 , 0.1,0.8), (y_5, 0.2,0.8) \}$\\
	$B_4=\{(y_1, 0.8,0.2), (y_2, 0.9,0.0), (y_3, 1.0,0.0), (y_4, 0.7,0.2), (y_5, 0.6,0.4) \}$
\end{table}
The symptoms of the given colorectal cancer patient $A$ are provided in terms of AIFS as follows:
\begin{table}[H]
	\scriptsize
	\centering
	$A=\{ (y_1, 0.3,0.5), (y_2, 0.4,0.4), (y_3, 0.6,0.2), (y_4, 0.5,0.5), (y_5, 0.9,0.0)\}$
\end{table}
\hspace*{-0.5cm}\textbf{Result Discussion:} The objective is to classify the state of the patient $A$ while using the known states $B_1,~ B_2$,~ $ B_3$, and $ B_4$. Let us use the proposed spatial similarity measure to calculate the similarities between the known states and $A$. The branches $\mbox{MD-}S_{\rm bv_{}}$, $\mbox{NMD-}S_{\rm bv_{}}$, and $\mbox{ED-}S_{\rm bv_{}}$ are classifying $A$ in the state $B_1$.
 The state of the patient $A$ is not identified by the similarity mesasures $S_{\rm L_1}, S_{\rm Hk}, S_{\rm M}, S^{2}_{\rm Hy_1}, S^{2}_{\rm Hy_2}, S^{2}_{\rm Hy_3}, S^{3}_{\rm Hy}, S_{\rm Ls}$ as their maximal similarity value is equal for two or three given patterns, whereas some existing similarity measures such as $S_{\rm L_2}, S_{\rm Lzd}, S_{\rm Bd}, S^{1}_{\rm Hy_1}, S^{1}_{\rm Hy_2}, S^{1}_{\rm Hy_3}, S_{\rm C}, S_{\rm Fz}, S_{\rm Dc}, S_{\rm Hm}, S_{\rm Az_{p}}, S_{\rm Az_{2}}^{h}$  recognize $A$ into state $B_1$ (see:\cite{2},\cite{23},\cite{24},\cite{55}). Thus, from Table-\ref{Example6}, we can say that our spatial similarity measure strongly classifies the patient $A$ to state $B_1$, i.e., Mestastasis.  \\
  \begin{table}[htp]
  	\centering
  	\fontsize{7}{7}\selectfont
  	\caption{Example 6: The classification obtained by the proposed SSM and other well-known measures for cancer diagnosis problem}
  		\renewcommand{\arraystretch}{1.5}\label{Example6}
  	\begin{tabular}{|p{2.5cm}|p{1.2cm}|p{1.2cm}|p{1.2cm}|p{1.2cm}|p{1.7cm}|}
  		\hline 
  		Similarity Measures &  $S(B_1,A)$ & $S(B_2,A)$ & $S(B_3,A)$&$S(B_4,A)$ &Result\\
  		\hline
  		$ S_{\rm L_{1}}$ & 0.8800 & 0.8800 & 0.5200 &0.6700 & Unclassified\\
  		\hline
  		$ S_{\rm Hk}$  & 0.8800 & 0.8800 & 0.5200 &0.6700 & Unclassified\\
  		\hline
  		$ S_{\rm L_{2}}$ & 0.8586 & 0.8388 & 0.4862 &0.6464 &$~~~~B_1$\\
  		\hline
  		$ S_{\rm Lzd}$  & 0.8586 & 0.8388 & 0.4862 &0.6464 &$~~~~B_1$\\
  		\hline
  		$ S_{\rm M}$~~($p=1$)  & 0.8800 & 0.8800 & 0.5200 &0.6700 &Unclassified\\
  		\hline
  		$ S_{\rm Dc}$~~($p=1$)& 0.9200 & 0.8800 & 0.5800 &0.6900 &$~~~~B_1$\\
  		\hline
  		$ S^{1}_{\rm Hy_1}$  & 0.8600 & 0.8200 & 0.4400 &0.6000 &$~~~~B_1$\\
  		\hline
 		$ S^{1}_{\rm Hy_2}$  &0.7933  &0.7394  &0.3217  &0.4785 &$~~~~B_1$\\
  		\hline
  		$ S^{1}_{\rm Hy_3}$  &0.7544  &0.6949  &0.2821  &0.4286 &$~~~~B_1$\\
  		\hline
  		$ S^{2}_{\rm Hy_1}$~~($p=1$) & 0.8800 & 0.8800 & 0.5200 &0.6700 &Unclassified\\
  		\hline
  		$ S^{2}_{\rm Hy_2}$~~($p=1$)  &0.7532  &0.7532  & 0.2863 &0.4412 &Unclassified\\
  		\hline
  		$ S^{2}_{\rm Hy_3}$~~($p=1$)&0.7097  & 0.7097 &0.2653  &0.4036 & Unclassified\\
  		\hline
  		$ S^{3}_{\rm Hy}$ & 0.8800 &0.8800  &0.5200  &0.6700 &Unclassified\\
  		\hline
  		$ S_{\rm Fz}$ &0.9000  &0.8800  &0.5500  &0.6800 &$~~~~B_1$\\
  		\hline
  		$ S_{\rm Bd}$~~($p=1,z=2$) &0.9067  &0.8800  & 0.5667 &0.6900 &$~~~~B_1$\\
  		\hline
  		$ S_{\rm C}$  &0.9200  &0.8800  &0.5800  &0.6900 &$~~~~B_1$\\
  		\hline
  		$ S_{\rm Ls}$~~($p=1$) & 0.8800 & 0.8800 & 0.5200 &0.6700 & Unclassified\\
  		\hline
  		$ S_{\rm Hm}$ & 0.8867 &   0.8600  &  0.5267 &   0.6533 &$~~~~B_1$\\
  		\hline
  		 $S_{\rm Az_{1}}$ &0.8625 &0.8125 &0.6375 &0.6875 &$~~~~B_1$ \\
  		\hline
  		$S_{\rm Az_{2}}$ &0.8380 &0.7331 &0.5655 &0.6779 &$~~~~B_1$ \\
  		\hline
  		$S_{\rm Az_{1}}^{h}$ &0.8750 &0.8000 &0.5500 &0.5125 &$~~~~B_1$ \\
  		\hline
  		$S_{\rm Az_{2}}^{h}$ &0.8459 &0.7764 &0.5417 &0.5191 &$~~~~B_1$ \\
  		\hline
  		 $\mbox{MD-}S_{\rm bv_{}}$ &0.9200 &0.8600 &0.8400 &0.8100 &$~~~~B_1$ \\
  		\hline
  		$\mbox{NMD-}S_{\rm bv_{}}$ &0.9350 &0.9150 &0.7650 &0.8150 &$~~~~B_1$ \\
  		\hline
  		$\mbox{ED-}S_{\rm bv_{}}$ &0.9450 &0.9250 &0.8550 &0.8750 &$~~~~B_1$ \\
  		\hline
  	\end{tabular}
  \end{table}

\section{Conclusion}\label{sec5}
The concept of bounded variation involves a difference operator. A novel divided difference operator of the double sequence is proposed to generalize the concept of bounded variation while introducing a multivalued absolutely summable bounded variation of two-dimensional AIFS. The multivalued absolutely summable bounded variation of two-dimensional AIFS is studied in terms of a proposed generalized metric space known as spatial metric space. The intrinsic adjustment of weights in spatial metric space is reflected in the form of three branches, i.e., membership dominant, non-membership dominant, and equidominant. In the two-dimensional AIFS domain, spatial metric space is a three branched distance measure that reduces to two branched distance measures in a fuzzy environment. The fuzzy set based two branches of distance measure is linearly dependent, hence only one branch helps in the pattern recognition. The spatial metric space reduces to the usual bounded variation based metric space in $\mathbb{R}.$ In the AIFS domain, the distance measure, and similarity measures complement each other; thus, we obtain a spatial similarity measure out of spatial metric space. The three branches of spatial similarity measure observe a maximum of three patterns simultaneously without involving weights. The weights are not required for SSM as they are intrinsically adjusted in the spatial metric space. The patterns recognized in the two-dimensional AIFS plane are more in comparison to the fuzzy or real domain, and it is reflected in some well-known pattern recognition problems, considered here.
The pattern recognition results obtained over numeric data problems, medical diagnosis, and cancer diagnosis with the proposed SSM are compared with well-known similarity measures to validate the performance of the measure. The complexity involved in classification while using the intuitionistic fuzzy similarity measures (see: Table \ref{similarity}) is the same as that of SSM and it is of order $n.$  

\section{Future Study}\label{sec6}
It is mainly a theoretical paper in which we have proposed a multivalued analysis in the two-dimensional AIFS domain. Here, a novel divided difference operator has been proposed in a two-dimensional AIFS domain. The proposed operator has been involved in the defining of spatial metric space. It is a fact that AIFS is a three-dimensional set, so the future requirement of this study is to propose a difference operator for AIFS. The proposal of a difference operator for AIFS will set the ground for the generalization of spatial metric space. The generalization of spatial metric space may give a direction to approach a well-known complex problem known as three body problem. A difference operator based sequence of bounded variation has been used to define the spatialistic view of the concept of bounded variation of sequences. A similar type of spatialistic study of other difference operator based sequence spaces may yield some interesting outcomes.\\
Spatial metric space has been utilized to classify the patterns in six benchmark pattern recognition problems. We suggest a rigorous application based study of SSM on real world problems to explore the inherent multiple patterns. The membership and non-membership components involving intuitionistic fuzzy similarity measures are considered for comparative purposes due to the proposal of multivalued analysis in the two-dimensional AIFS plane. In the future, we intend to propose spatial decision making methods independent of criteria weights.\\

\hspace*{-0.5cm}\textbf{Disclosure of Competing Interests}\\
The authors verify that they do not have any known competing financial interests or personal relationships that could have potentially impacted the findings presented in this paper.\\

\hspace*{-0.5cm}\textbf{Details of Funding}\\
This research received no specific grants from anywhere.

\begin{spacing}{0.1}

\end{spacing}

\end{document}